\documentclass{amsart}

\usepackage{amsmath,amsthm,amsfonts,amscd,amssymb,mathabx}
\usepackage[pdftex]{graphicx}
\usepackage{array}
\usepackage[numbers]{natbib}
\usepackage{verbatim}
\usepackage{enumerate}
\usepackage{fullpage}
\usepackage[dvipsnames]{xcolor}
\usepackage[pdfborder={0 0 0},colorlinks,linkcolor=blue,citecolor=magenta,urlcolor=black,hypertexnames=false]{hyperref} 

\newcommand{\ZZ}{\mathbb{Z}}

\newtheorem{thm}{Theorem}[section]
\newtheorem{cor}[thm]{Corollary}
\newtheorem{lem}[thm]{Lemma}
\newtheorem{prop}[thm]{Proposition}

\newtheorem{remark}{Remark}[section]

\theoremstyle{remark}
\newtheorem{ex}{Example}

\theoremstyle{definition}

\theoremstyle{remark}
\newtheorem{rem}{Remark}[section]

\newenvironment{psmatrix}
  {\left(\begin{smallmatrix}}
  {\end{smallmatrix}\right)}

\newcolumntype{P}[1]{>{\centering\arraybackslash}p{#1}}

\begin{document}

\title[Lattices From Tight Frames and Vertex Transitive Graphs]{Lattices From Tight Frames and Vertex Transitive Graphs}
\author{Lenny Fukshansky}\thanks{LF was partially supported by the Simons Foundation grant \#519058}
\author{Deanna Needell}\thanks{DN was partially supported by the NSF CAREER DMS \#1348721 and the NSF BIGDATA DMS \#1740325}
\author{Josiah Park}\thanks{JP was partially supported by grant from the NSF DMS \#1600693}
\author{Yuxin Xin}

\address{Department of Mathematics, 850 Columbia Avenue, Claremont McKenna College, Claremont, CA 91711}
\email{lenny@cmc.edu}
\address{Department of Mathematics, University of California at Los Angeles, 520 Portola Plaza, Los Angeles, CA 90095}
\email{deanna@math.ucla.edu}
\address{School of Mathematics, Georgia Institute of Technology, 686 Cherry Street, Atlanta, GA 30332-0160}
\email{j.park@gatech.edu}
\address{Department of Mathematics, 850 Columbia Avenue, Claremont McKenna College, Claremont, CA 91711}
\email{yxin19@students.claremontmckenna.edu}


\begin{abstract}
We show that real tight frames that generate lattices must be rational, and use this observation to describe a construction of lattices from vertex transitive graphs. In the case of irreducible group frames, we show that the corresponding lattice is always strongly eutactic. This is the case for the more restrictive class of distance transitive graphs. We show that such lattices exist in arbitrarily large dimensions and demonstrate examples arising from some notable families of graphs. In particular, some well-known root lattices and those related to them can be recovered this way. We discuss various properties of this construction and also mention some potential applications of lattices generated by incoherent systems of vectors.
\end{abstract}

\maketitle

\def\A{{\mathcal A}}
\def\AA{{\mathfrak A}}
\def\B{{\mathcal B}}
\def\C{{\mathcal C}}
\def\D{{\mathcal D}}
\def\EE{{\mathfrak E}}
\def\F{{\mathcal F}}
\def\G{{\mathcal G}}
\def\x{{\mathcal H}}
\def\I{{\mathcal I}}
\def\II{{\mathfrak I}}
\def\J{{\mathcal J}}
\def\K{{\mathcal K}}
\def\kk{{\mathfrak K}}
\def\L{{\mathcal L}}
\def\LL{{\mathfrak L}}
\def\M{{\mathcal M}}
\def\mm{{\mathfrak m}}
\def\MM{{\mathfrak M}}
\def\N{{\mathcal N}}
\def\O{{\mathcal O}}
\def\OO{{\mathfrak O}}
\def\PP{{\mathfrak P}}
\def\R{{\mathcal R}}
\def\W{{\mathcal W}}
\def\PNR{{\mathcal P_N(\real)}}
\def\PMNR{{\mathcal P^M_N(\real)}}
\def\PdNR{{\mathcal P^d_N(\real)}}
\def\s{{\mathcal S}}
\def\V{{\mathcal V}}
\def\X{{\mathcal X}}
\def\Y{{\mathcal Y}}
\def\Z{{\mathcal Z}}
\def\Z{{\mathcal Z}}
\def\({{\langle}}
\def\){{\rangle}}
\def\cee{{\mathbb C}}
\def\Nn{{\mathbb N}}
\def\pee{{\mathbb P}}
\def\que{{\mathbb Q}}
\def\QQ{{\mathbb Q}}
\def\real{{\mathbb R}}
\def\RR{{\mathbb R}}
\def\zed{{\mathbb Z}}
\def\ZZ{{\mathbb Z}}
\def\aaa{{\mathbb A}}
\def\ff{{\mathbb F}}
\def\HDelta{{\it \Delta}}
\def\kk{{\mathfrak K}}
\def\qbar{{\overline{\mathbb Q}}}
\def\kbar{{\overline{K}}}
\def\ybar{{\overline{Y}}}
\def\kkbar{{\overline{\mathfrak K}}}
\def\ubar{{\overline{U}}}
\def\eps{{\varepsilon}}
\def\ahat{{\hat \alpha}}
\def\bhat{{\hat \beta}}
\def\k{{\nu}}
\def\gt{{\tilde \gamma}}
\def\h{{\tfrac12}}
\def\be{{\boldsymbol e}}
\def\bei{{\boldsymbol e_i}}
\def\bc{{\boldsymbol c}}
\def\bdt{{\boldsymbol \delta}}
\def\bff{{\boldsymbol f}}
\def\bm{{\boldsymbol m}}
\def\bk{{\boldsymbol k}}
\def\bh{{\boldsymbol h}}
\def\bi{{\boldsymbol i}}
\def\bl{{\boldsymbol l}}
\def\bp{{\boldsymbol p}}
\def\bq{{\boldsymbol q}}
\def\br{{\boldsymbol r}}
\def\bu{{\boldsymbol u}}
\def\bt{{\boldsymbol t}}
\def\bs{{\boldsymbol s}}
\def\bv{{\boldsymbol v}}
\def\bw{{\boldsymbol w}}
\def\bx{{\boldsymbol x}}
\def\bb{{\boldsymbol b}}
\def\bbx{{\overline{\boldsymbol x}}}
\def\bX{{\boldsymbol X}}
\def\bz{{\boldsymbol z}}
\def\bwy{{\boldsymbol y}}
\def\bwv{{\boldsymbol v}}
\def\bY{{\boldsymbol Y}}
\def\bL{{\boldsymbol L}}
\def\ba{{\boldsymbol a}}
\def\bb{{\boldsymbol b}}
\def\bet{{\boldsymbol\eta}}
\def\bxi{{\boldsymbol\xi}}
\def\bo{{\boldsymbol 0}}
\def\bone{{\boldsymbol 1}}
\def\bol{{\boldsymbol 1}_L}
\def\ep{\varepsilon}
\def\p{\boldsymbol\varphi}
\def\q{\boldsymbol\psi}
\def\rank{\operatorname{rank}}
\def\aut{\operatorname{Aut}}
\def\lcm{\operatorname{lcm}}
\def\sgn{\operatorname{sgn}}
\def\spn{\operatorname{span}}
\def\md{\operatorname{mod}}
\def\Norm{\operatorname{Norm}}
\def\dim{\operatorname{dim}}
\def\det{\operatorname{det}}
\def\Vol{\operatorname{Vol}}
\def\rk{\operatorname{rk}}
\def\ord{\operatorname{ord}}
\def\ker{\operatorname{ker}}
\def\div{\operatorname{div}}
\def\Gal{\operatorname{Gal}}
\def\GL{\operatorname{GL}}
\def\SNR{\operatorname{SNR}}
\def\WR{\operatorname{WR}}
\def\IWR{\operatorname{IWR}}
\def\scg{\operatorname{\left< \Gamma \right>}}
\def\swrh{\operatorname{Sim_{WR}(\Lambda_h)}}
\def\ch{\operatorname{C_h}}
\def\cht{\operatorname{C_h(\theta)}}
\def\scgt{\operatorname{\left< \Gamma_{\theta} \right>}}
\def\scgmn{\operatorname{\left< \Gamma_{m,n} \right>}}
\def\gat{\operatorname{\Omega_{\theta}}}
\def\Obar{\operatorname{\overline{\Omega}}}
\def\Lbar{\operatorname{\overline{\Lambda}}}
\def\mn{\operatorname{mn}}
\def\disc{\operatorname{disc}}
\def\rot{\operatorname{rot}}
\def\Prob{\operatorname{Prob}}
\def\co{\operatorname{co}}
\def\ot{\operatorname{o_{\tau}}}
\def\Aut{\operatorname{Aut}}
\def\Mat{\operatorname{Mat}}
\def\SL{\operatorname{SL}}
\def\id{\operatorname{id}}
\def\skel{\operatorname{skel}}

\section{Introduction}
\label{intro}

Let $\left<\ ,\ \right>$ be the usual inner product on $\real^k$ and $\|\bx\| := \left< \bx,\bx \right>^{1/2}$ the Euclidean norm on $\mathbb{R}^k$.
For a lattice $L \subset \real^k$ of full rank $k$ (that is a discrete co-compact subgroup of $\real^k$) the {\it minimal norm} of $L$ is
$$|L| := \min \{ \|\bx\| : \bx \in L \setminus \{ \bo \} \},$$
and its set of {\it minimal} or {\it shortest vectors} is
$$S(L) := \{ \bx \in L : \|\bx\| = |L| \}.$$
The {\it automorphism group} of the lattice $L$, $\Aut(L)$, is the group of all $k \times k$ real orthogonal matrices that map $L$ to itself. A particularly interesting class of lattices are {\it eutactic} lattices: a lattice $L$ is called eutactic if its set of minimal vectors $S(L)$ satisfies a eutaxy condition, i.e. there exist positive real numbers $c_1,\dots,c_n$, (called eutaxy coefficients) such that
\begin{equation}
\label{eutaxy}
\|\bv\|^2 = \sum_{\bx \in S(L)} c_i \left< \bv,\bx_i \right>^2
\end{equation}
for all $\bv \in \real^k$. If $c_1 = \dots = c_n$, $L$ is said to be {\it strongly eutactic}. Eutactic and strongly eutactic lattices are central objects of lattice theory due to their importance in connection with well studied optimization problems. A theorem of Voronoi (1908) asserts that $L$ is a local maximum of the packing density function on the space of lattices in $\real^k$ if and only if $L$ is eutactic and perfect ($L$ is {\it perfect} if the set $\{ \bx^{\top} \bx : \bx \in S(L) \}$ spans the space of $k \times k$ real symmetric matrices)~\cite{voronoi}. More details on eutactic, strongly eutactic and perfect lattices can be found in J. Martinet's book~\cite{martinet}.

Two lattices $L$ and $M$ are called {\it similar}, written $L \sim M$, if $L = \alpha U M$ for a nonzero scalar $\alpha$ and an orthogonal transformation $U$. Similarity is an equivalence relation on lattices that preserves inner products between vectors (up to the scalar $\alpha$) and, as a result, lattice's automorphism group; it also gives a bijection between sets of minimal vectors. Consequently, all the geometric properties that we discuss here, such as eutaxy, strong eutaxy and perfection are preserved on similarity classes.

In the previous papers~\cite{etf1} and~\cite{etf2} of the first two authors, lattices generated by equiangular tight frames (ETFs) were studied and examples of strongly eutactic such lattices were constructed. Here we aim to take this discussion further. Let $n \ge k$ and let $\F := \left\{ \bff_1,\dots,\bff_n \right\} \subset \real^k$ be a sequence of vectors, not necessarily distinct, such that
$\spn_{\real} \left\{ \bff_1,\dots,\bff_n \right\} =\real^k$.
Such a set $\F$ is called an $(n,k)$-{\it frame}, the name originating in a 1952 paper of Duffin and Schaeffer in connection with their study of nonharmonic Fourier series~\cite{duffin}. A frame $\F$ is called {\it uniform} if all of its vectors have the same norm, and it is called {\it tight} if there exists a real constant $\gamma > 0$ such that for every $\bv \in \real^k$
\begin{equation}
\label{tight}
\|\bv\|^2 = \gamma \sum_{i=1}^n \langle \bv,\bff_i \rangle^2,
\end{equation}
and a tight frame is called {\it Parseval} if $\gamma=1$: clearly, any tight frame can be rescaled to a Parseval frame. Notice the similarity between this equation and the equation~\eqref{eutaxy} above. Although the tightness condition~\eqref{tight} above is well studied in several contemporary branches of mathematics, the closely related eutaxy condition precedes it by half a century. Voronoi's study~\cite{voronoi} of quadratic forms in 1908 gave rise to the introduction of eutaxy condition~\eqref{eutaxy}. Nonetheless, we can say that a lattice is strongly eutactic whenever its set of minimal vectors forms a uniform tight frame. Another way to view uniform tight frames is as projective $1$-designs, a subclass of more general designs on compact spaces introduced by Delsarte, Goethals, and Seidel in their groundbreaking 1977 paper~\cite{delsarte}. A special class of tight frames are examples of optimal packings of lines in projective space. These uniform tight frames are called {\it equiangular} (abbreviated ETF) if $\left| \left< \bff_i, \bff_j \right> \right|$ is the same for all $i \neq j$. Tight frames in general and ETFs in particular are extensively studied objects in harmonic analysis; see S. Waldron's book~\cite{waldron_book} for detailed information on this subject.

Given a real $(n,k)$-frame $\F = \left\{ \bff_1,\dots,\bff_n \right\}$, define
$$L(\F) = \spn_{\zed} \left\{ \bff_1,\dots,\bff_n \right\}.$$
If we write $B$ for the $k \times n$ matrix with vectors  $\bff_1,\dots,\bff_n$ as columns, then
$$L(\F) = \{B \ba: \ba \in \zed^n\}.$$
The {\it norm-form} associated with $\F$ is the quadratic form
\begin{equation}
\label{QF}
Q_\F(\ba)=\| B \ba\|^2=\langle B^\top B \ba, \ba \rangle.
\end{equation}
We call the frame $\F$ {\it rational} if $Q_\F$ is (a constant multiple of) a rational quadratic form, i.e.  the $n \times n$ Gram matrix $B^\top B$ is (a constant multiple of) a rational matrix. This is equivalent to saying that the inner products $\left< \bff_i, \bff_j \right>$ are (up to a constant multiple) rational numbers for all $1 \leq i,j \leq n$. In~\cite{etf2}, it was proved that if $\F$ is rational, then $L(\F)$ is a lattice. Further, in the case that $\F$ is an ETF, $L(\F)$ is a lattice if and only if $\F$ is rational (the converse was previously proved in~\cite{etf1}). More generally, it was shown in~\cite{etf2} that when the dimension $k=2$ or $3$ and $\F$ is a tight $(n,k)$-frame for any $n$ so that $L(\F)$ is a lattice, then $\F$ must be rational. Our first result is an extension of this observation to any dimension.

\begin{thm} \label{rational} Suppose that $\F$ is a tight $(n,k)$-frame so that $L(\F)$ is a lattice. Then $\F$ must be rational.
\end{thm}

\noindent
We give two different proofs of Theorem~\ref{rational} in Section~\ref{rat}, one of them as a consequence of a stronger result about a larger class of matrices than just the tight frames (Theorem~\ref{thm_rat}).

\begin{rem} \label{sunada_rem} We have recently become aware of a 2017 paper by T. Sunada~\cite{sunada}, where a result similar to our Theorem~\ref{rational} has been established (Proposition~4.2 of~\cite{sunada}). This being said, our proof of this result is considerably simpler, and our Theorem~\ref{thm_rat} is more general: it does not follow from~\cite{sunada}.
\end{rem}

We can now use this rationality result to pick out lattices generated by tight frames. We are especially interested in frames that give rise to lattices with nice geometric properties. For this we need some more notation. Let the automorphism group of a frame $\F$ be
$$\Aut(\F) := \{ U \in \O_k(\real) : U \bff \in \F \text{ for all } \bff \in \F \},$$
where $\O_k(\real)$ is the group of $k \times k$ real orthogonal matrices. As usual, we write $H \leq G$ to indicate that $H$ is a subgroup of the group~$G$.

We now discuss {\it group frames}; see Chapter~10 of~\cite{waldron_book} for a detailed exposition. Let $\bff_1 \in \real^k$ be a vector and let $G$ a finite  group of orthogonal $k \times k$ matrices. Define $\F$ to be the orbit of $\bff_1$ under the action of $G$ by left multiplication, i.e.
$$\F = G\bff_1:= \left\{ U \bff_1 : U \in G \right\},$$
then all the vectors in $\F$ have the same norm. If $\F$ spans $\real^k$, then $\F$ is a uniform frame, which we refer to as a $G$-{\it frame}. $G$ is said to act {\it irreducibly} on the space~$\real^k$ if there is no nonzero proper subspace $E$ of $\real^k$ that is closed under the action of $G$, that is, $GE \neq E$ for any $\{ \bo \} \neq E \subsetneq \real^k$. A $G$-frame with such an irreducible action corresponding to $G$ on~$\real^k$ is similarly called irreducible. All irreducible group frames are tight. In fact, if $G$ is a group with an irreducible action on~$\real^k$, then the orbit of $x$ under $G$, $\{ U \bx : U \in G \}$, is an irreducible tight $G$-frame for any nonzero vector $\bx \in \real^k$ (see Sections~10.5 - 10.9 of~\cite{waldron_book} for details). 

Our next result demonstrates a certain correspondence between irreducible group frames and strongly eutactic lattices.

\begin{thm} \label{Gframes} Let $G$ be a group of $k \times k$ real orthogonal matrices and $\bff \in \real^k$ be a vector so that $\F = G\bff$ is an irreducible rational group frame in~$\real^k$. Then the lattice $L(\F)$ is strongly eutactic. 
\end{thm}

\begin{rem} \label{thm2_converse} Conversely, suppose $L \subset \real^k$ is a strongly eutactic lattice of rank~$k$. By Corollary~16.1.3 of~\cite{martinet}, $L$ is strongly eutactic if and only if its set $S(L)$ of minimal vectors is a spherical 2-design, which is a condition equivalent to the tightness condition~\eqref{tight}. Since all minimal vectors have the same norm, $S(L)$ is a uniform tight frame. Now suppose some $\Aut(L)$ acts transitively on $S(L)$. Let $\bx_1 \in S(L)$, then for any $\bx \in S(L)$ there exists a $U \in \Aut(L)$ such that $\bx = U \bx_1$. Hence
$$S(L) = \{ U \bx_1 : U \in \Aut(L) \},$$
and so $S(L)$ is an $\Aut(L)$-frame. If the action of~$\Aut(L)$ on~$\real^k$ is irreducible then $S(L)$ is an irreducible group frame.
\end{rem}
We prove Theorem~\ref{Gframes} in Section~\ref{grp_frames}. This theorem motivates the investigation of rational irreducible group frames. One steady source of rational group frames comes from vertex transitive graphs, as detailed in Section~10.7 of~\cite{waldron_book}. In the special case when the graph in question is distance transitive, these frames are irreducible.

\begin{thm} \label{graph_thm}
Let $\Gamma$ be a vertex transitive graph on $n$ vertices and $G$ its automorphism group. Let $A$ be the adjacency matrix of $\Gamma$ and $\lambda$ an eigenvalue of multiplicity~$m$. Assume $\lambda$ is rational and let $V_{\lambda}$ be the corresponding $m$-dimensional eigenspace to eigenvalue $\lambda$. Let $P_{\lambda}$ be a rational orthogonal projection matrix of~$\real^n$ onto~$V_{\lambda}$. Then $L_{\Gamma,\lambda} := P_{\lambda} \zed^n$ is a lattice of full rank in $V_{\lambda}$, and its automorphism group contains a subgroup isomorphic to a factor group of~$G$. If $\Gamma$ is distance transitive, $L_{\Gamma,\lambda}$ is strongly eutactic.
\end{thm}

We review all the necessary notation and prove Theorem~\ref{graph_thm} in Section~\ref{graphs}. Distance transitive graphs form a subclass of vertex transitive graphs, and there are plenty of examples of such graphs with rational eigenvalues. In fact, there exist such lattices on $n$ vertices for arbitrarily large $n$ having eigenvalues of multiplicity $m$ being an increasing function of $n$ (for instance complete graphs, Johnson graphs, Grassman graphs, folded cube graphs, etc.), so that this construction yields strongly eutactic lattices in arbitrarily high dimensions. Further, there are some instances of vertex transitive graphs which are not distance transitive, however still give rise to strongly eutactic lattices. We demonstrate several examples of our construction in Section~\ref{graphs}, some of which are summarized in Table~\ref{graph_ex}. A separate collection of lattices coming from several Johnson graphs $J(n,2)$ is given in Table~\ref{J_graph_ex} in Section~\ref{graphs}. Furthermore, in Theorem~\ref{tensor} we give a characterization of lattices coming from product graphs in terms of tensor products and orthogonal direct sums of component lattices.

For the purposes of all of our examples and constructions, the lattices are viewed up to similarity and eigenspaces of graphs are identified with real Euclidean spaces~$\real^k$ for the appropriate dimension $k$ equal to the multiplicity of the corresponding eigenvalue. Our examples have been computed in Maple~\cite{maple} using online catalog~\cite{vt_graph} of distance regular graphs and online catalog~\cite{streut} of strongly eutactic lattices. It can be seen from these examples that a graph and its complement produce the same lattices. This is true in general, as is shown in Proposition~\ref{complement} in Section~\ref{graphs}. At the end of Section~\ref{graphs} we also demonstrate an interesting correspondence between contact polytopes of lattices $E_6^*$, $E_7^*$ and $A_3^*$ and our construction of lattices from their skeleton graphs.

\begin{center}
\begin{table} 
\begin{tabular}{|P{4.5cm}|P{1.5cm}|P{1.5cm}|P{1.5cm}|P{1.5cm}|P{3.5cm}|} \hline
{\bf Graph $\Gamma$} & {\bf Dist. trans.?} & {\bf \# of vertices} & {\bf Eig. $\lambda$} & {\bf Mult. of $\lambda$} & {\bf Lattice $L_{\Gamma,\lambda}$} \\
\hline \hline
Disconnected graph & No & ($n$) & $\ \ 0$ & ($n$) & Integer lattice $\zed^n$ \\ \hline
Complete graph $K_n$ & Yes & ($n$) & $-1$ & ($n-1$) & Root lattice $A_{n-1}$ \\ \hline
Hamming graph $H(2,3)$ & Yes & ($9$) & $\ \ 1$ & ($4$) & $A_2 \otimes_{\zed} A_2$ \\ \hline
Petersen graph & Yes & ($10$) & $-2$ & ($4$) & $A_4^*$, dual of $A_4$ \\ \hline
Petersen graph & Yes & ($10$) & $\ \ 1$ & ($5$) & Coxeter lattice $A_5^2$ \\ \hline
Petersen line graph & Yes & ($15$) & $-1$ & ($4$) & $A_4^*$, dual of $A_4$ \\ \hline
Petersen line graph & Yes & ($15$) & $-2$ & ($5$) & Coxeter lattice $A_5^3$ \\ \hline
Clebsch graph & Yes & ($16$) & $-3$ & ($5$) & $D_5^*$, dual of $D_5$ \\ \hline
Clebsch graph complement & Yes & ($16$) & $\ \ 2$ & ($5$) & $D_5^*$, dual of $D_5$ \\ \hline
Shrikhande graph & No & ($16$) & $\ \ 2$ & ($6$) & $D_6^+$ \\ \hline
Shrikhande complement & No & ($16$) & $-3$ & ($6$) & $D_6^+$ \\ \hline
Schl\"afli graph & Yes & ($27$) & $\ \ 4$ & ($6$) & $E_6^*$, dual of $E_6$ \\ \hline
Schl\"afli graph complement & Yes & ($27$) & $-5$ & ($6$) & $E_6^*$, dual of $E_6$ \\ \hline
Gosset graph & Yes & ($56$) & $\ \ 9$ & ($7$) & $E_7^*$, dual of $E_7$ \\ \hline
\end{tabular}
\medskip
\caption{Examples of strongly eutactic lattices from vertex transitive graphs} 
\label{graph_ex}
\end{table}
\end{center}

It is also interesting to consider Theorem~\ref{graph_thm} in view of the properties of eutactic configurations, i.e. finite sets of vectors satisfying the eutaxy condition~\eqref{eutaxy}. The famous theorem of Hadwiger (\cite{martinet}, Theorem~3.6.12) asserts that a set $S$ of cardinality $n$ in $k$-dimensional space $V$, $n > k$, is eutactic if and only if it is an orthogonal projection onto $V$ of an orthonormal basis in an $n$-dimensional space containing~$V$. In fact, our construction considers precisely such a projection, namely the set of vectors $\{ P_{\lambda} \be_i \}_{i=1}^n$ where $\be_1,\dots,\be_n$ is the standard basis in~$\real^n$. This set is therefore eutactic by Hadwiger. Our result, however, implies more, specifically that in our setting (in the case of distance transitive graphs) these vectors generate a lattice whose set of minimal vectors is strongly eutactic.

Finally, in Section~\ref{coh} we discuss a possible relation between coherence of a lattice and its sphere packing density, as well as potential applications of tight frames coming from sets of minimal vectors of lattices in compressed sensing.
\bigskip

\section{Rationality of lattice-generating frames}
\label{rat}

We start with a simple proof of Theorem~\ref{rational}.

\proof[Proof of Theorem~\ref{rational}]
With notation as in the statement of the theorem, let $B$ be a $k \times n$ real matrix whose columns are vectors of the tight frame~$\F$ and $L(\F)$ is a lattice. Let $A$ be a $k \times k$ basis matrix for $L(\F)$. Then, there exists a $k \times n$ integer matrix $Z$ so that $AZ = B$. Thus 
$$AZZ^{\top}A^{\top} = BB^{\top} = \gamma I_k$$
for some $\gamma>0$. Since $A$ is invertible, 
$$ZZ^{\top} = \gamma A^{-1}(A^{\top})^{-1},$$
so that $ZZ^{\top}=\gamma (A^{\top}A)^{-1}$. Therefore 
$$B^{\top} B = Z^{\top} A^{\top} A Z = Z^{\top} \gamma (ZZ^{\top})^{-1} Z = \gamma Z^{\top} (ZZ^{\top})^{-1} Z.$$
Since $Z^{\top}(ZZ^{\top})^{-1} Z$ has rational entries, we have that $B^{\top} B$ is a multiple of a rational matrix. Therefore $\F$ is a rational tight frame.
\endproof

The above argument implies that if $Q_{\F}$ as in~\eqref{QF} is a quadratic form corresponding to an irrational tight frame~$\F$ then the corresponding integer span~$L(\F)$ is not a lattice (i.e. is not discrete) because~$Q_{\F}$ cannot be bounded away from $0$ on integer points. This argument, however, relies heavily on the norm-form $Q_{\F}$ coming from a tight frame. On the other hand, it is not difficult to construct other irrational quadratic forms (not corresponding to tight frames) which are bounded away from $0$ on integer points. For instance, take $L_1,\dots,L_k$ to be rational linear forms in $n$ variables $x_1,\dots,x_n$ and $c_1,\dots,c_k$ any positive real numbers. Let
$$Q(x_1,...,x_n) = c_1 L_1^2 + \dots + c_k L_k^2.$$
This $Q$ is a positive semidefinite quadratic form. Suppose $Q(\ba) \neq 0$ for some integer vector~$\ba$, then there must exist $1 \leq i \leq k$ such that $L_i(a) \neq 0$. Since $L_i$ has rational coefficients, $|L_i(a)| \geq 1/d_i$, where $d_i$ is the least common multiple of the denominators of these coefficients. Let $d = \max \{ d_1,\dots,d_k \}$ and $c = \min \{ c_1,\dots,c_k \}$, then we have
$$Q(a) \geq c/d^2$$
for all $\ba$ for which $Q(\ba) \neq 0$. In particular, if some of the $c_i$'s are irrational, $Q$ is a form with irrational coefficients.

In view of this observation, it is interesting to understand what are the necessary and sufficient conditions on a $k \times n$ real matrix $B$ so that $B \zed^n$ is a lattice to imply that $B$ must be rational? In the rest of this section we prove a sufficient condition that is weaker than being a tight frame. Write $\{\bb_i\}_{i=1}^n \subset \real^k$ for the elements of a frame $\F$ (a sequence of vectors spanning $\real^k$), written as column vectors of a $k \times n$ matrix $B$, where $n=k+m$. Let the first $k$ columns in $B$ be denoted in matrix form by $B_{0}$ and the remaining $m$ column vectors by $B_{1}$, so that $B=[B_{0}\ |\ B_{1}]$, $B_{0} \in \real^{k \times k}$,\ $B_{1} \in \real^{k \times m}$. 

\begin{lem} \label{invrat} Suppose that $B=[B_{0}\ |\ B_{1}]$ is such that $B_{0} \real^k = \real^k$ and $\Lambda_{B} := B \zed^n$ is discrete. Then $B_{0}^{-1}B_{1}\in \que^{k \times m}$.
\end{lem}

\proof
If $\Lambda_{B}$ is discrete, it is a full-rank lattice in~$\real^k$, and so has a basis matrix $A = \begin{pmatrix} \ba_1 & \dots & \ba_k \end{pmatrix}$ such that $\Lambda_{B} = A\zed^{k}$. Hence there exist some integer matrices $Z_{0},Z_{1}$ such that $AZ_{0}=B_{0}$, and $AZ_{1}=B_{1}$. Since $B_{0}$ is full rank and $A$ invertible, $Z_{0}$ is invertible and $B_{0}^{-1}B_{1}=Z_{0}^{-1}A^{-1}AZ_{1}=Z_{0}^{-1}Z_1\in \que^{k \times m}$.
\endproof

Let $Q$ be an $k \times k$ orthogonal real matrix, then multiplication by $Q$ preserves inner products of vectors in~$\real^k$ and a collection of vectors $\{\bb_i\}_{i=1}^n$ generates a lattice over $\zed$ if and only if $\{Q \bb_i\}_{i=1}^n$ does. Let $W$ be orthogonally equivalent to $B$, that is $W=QB$ for some $Q \in \O_k(\real)$ ($\O_k(\real)$ denotes the set of real $k\times k$ orthogonal matrices). $QQ^{\top} = I_k$, the $k \times k$ identity matrix, and the matrix of outer products for $W$ is $WW^{\top} = QBB^{\top} Q^{\top}$. Having information about the entries of this matrix for certain $Q$ (arising in this case from the $QR$-decomposition of a matrix) allows for an easy way to check rationality of inner products. When $B$ is a tight frame given in matrix form, (as above) $BB^{\top} = \gamma I_k$ for some $\gamma>0$, and so $WW^{\top}$ collapses to the same matrix as $BB^{\top}$.  In general, however the relationship between $WW^{\top}$ and $BB^{\top}$ can get ``muddled" by transformation so that determining lattice properties of integer combinations of vectors in a tight frame is easier than the general case.

\begin{remark} Given $B_{0}=QR$, the $QR$ factorization of $B_{0}$, so that $Q \in \O_k(\real)$ and $R$ is upper-triangular with positive entries along the diagonal, it will be useful to work with the alternative representation of $B$: $\tilde{B}=Q^{-1}B=[R\ |\ Q^{-1}B_{1}]$. In the arguments which follow, we choose to write $\tilde{B}=D[U\ |\ V]$, where $D\in\real^{k\times k}$ is diagonal with entries $d_1,...,d_k$, $U\in \real^{k\times k}$ is upper-triangular with ones along the diagonal, and $V\in\real^{k\times m}$ is the remaining entries. In the above, $d_i$ are taken to be positive (which is possible since $R$ has positive diagonal entries). From now on, let $B$ denote a matrix of the form $\tilde{B}$ when not specified otherwise. 
\end{remark}

\begin{thm} \label{thm_rat} Suppose a collection of vectors $B=\{ \bb_{i}\}_{i=1}^{n} \subset \real^{k}$, $n=k+m$, is given as column vectors of a matrix of the form $\tilde{B}$ (as in the preceding remark). Suppose these column vectors have the following properties:
 \begin{enumerate}[(i)]

\item $\spn_{\zed} B$ is discrete, 

\item the row-vectors of $B$, $\br_1,\dots,\br_k$, satisfy $\langle \br_{i}, \br_{j} \rangle = d_{i}d_{j}q_{i,j}$ for some $q_{i,j}\in \que$ and all $i\neq j$, that is, $[U\ |\ V][U\ |\ V]^{T}$ has rational entries off the diagonal, and

\item $\langle \br_{i}, \br_{i} \rangle=q_{i,i}\in\que$ for all $i=1,\dots,k$, that is, $BB^{T}=D[U\ |\ V][U\ |\ V]^{T}D$ has rational entries on the diagonal.

\end{enumerate}
Then the inner products $\langle \bb_{i}, \bb_{j} \rangle$ must all be rational, i.e. $B^{\top} B \in \que^{n \times n}$.
\end{thm}

\proof
For each column vector $\bb_{j}$ from $B$, Lemma~\ref{invrat} implies there exists a vector $\bp_{j}\in\que^{k}$ such that $B_{0}^{-1} \bb_{j}= \bp_{j}$. Letting $p_{i,j}$ be the $i$-th entry of each $\bp_{j}$, we now use these rational numbers to demonstrate, under the above conditions, that $B^{\top} B$ must be rational.

\begin{figure} $ \protect\begin{pmatrix} d_{1} & d_{1}u_{1,2} & d_{1}u_{1,3} & d_{1}u_{1,4} & \hdots & d_{1}v_{1,1} & d_{1}v_{1,2} & \hdots & d_{1}v_{1,m} \\ 
0 & d_2 & d_2u_{2,3} & d_2u_{2,4} & \hdots & d_2v_{2,1} & d_2v_{2,2} & \hdots & d_2v_{2,m} \\ 
0 & 0 & d_3 & d_3u_{3,4} & \hdots & d_3v_{3,1} & d_3v_{3,2} & \hdots & d_3v_{3,m} \\ 
\vdots & \vdots & \vdots & \vdots & \ddots & \vdots & \vdots & \ddots & \vdots \\ 
0 & 0 & 0 & 0 & \hdots & d_kv_{k,1} & d_kv_{k,2} & \hdots & d_kv_{k,m} \\    \protect\end{pmatrix} $ \caption{Matrix $B$.} \end{figure}

Recall that $B$ has $k$ rows and $k+m$ columns. From now on, denote the last $m$ column vectors of $B$ by $\bwv_{l}$, $1\leq l \leq m$. The condition $B_{0}^{-1} \bwv_{l} = \bp_{l}$ gives that $d_k p_{k,l} = d_kv_{k,l}$, so $v_{k,l} = p_{k,l}$ for all $l = 1,\dots,m$, i.e., the numbers $v_{k,l}$ are rational. In the same manner, we obtain $m$ equations:
$$d_{k-1} p_{k-1,l}+d_{k-1}u_{k-1,k}p_{k,l} = d_{k-1}v_{k-1,l},$$
which imply that $p_{k-1,l}+u_{k-1,k}p_{k,l} = v_{k-1,l}$ for all $l = 1,\dots,m$, as well as 
$$d_kd_{k-1}(u_{k-1,k}+v_{k-1,1}v_{k,1}+\dots+v_{k-1,m}v_{k,m}) = q_{k,k-1}d_kd_{k-1},$$
which implies $u_{k-1,k}+v_{k-1,1}v_{k,1}+\dots+v_{k-1,m}v_{k,m} = q_{k,k-1}$. Now, these $m+1$ equations can be written together in a matrix equation:
$$\begin{pmatrix}   1 & v_{k,1} & v_{k,2} & \hdots & v_{k,m} \\ 
p_{k,1} & -1 & 0 & \hdots & 0 \\
p_{k,2} & 0 & -1 & \hdots & 0 \\
\vdots  & \vdots & \vdots & \ddots & \vdots \\
p_{k,m} & 0 & 0 & \hdots & -1 \\
 \end{pmatrix}\begin{pmatrix} u_{k-1,k} \\ v_{k-1,1} \\ v_{k-1,2} \\ \vdots \\ v_{k-1,m} \end{pmatrix} = \begin{pmatrix} q_{k,k-1} \\ -p_{k-1,1} \\ -p_{k-1,2} \\ \vdots \\ -p_{k-1,m} \end{pmatrix}.$$
The matrix formed above on the left is invertible as all the rows of index greater than one are orthogonal to the first (this may be checked using the condition $v_{k,l}=p_{k,l}$) and the lower right block being the negative identity shows the last $m$ rows (arising from the first set of equalities above) are linearly independent amongst themselves. Thus by applying the inverse of the matrix on the left to each side we can express the coordinates $u_{k-1,k},v_{k-1,1},\dots,v_{k-1,m}$ of the vector on the left as rational numbers.

Proceeding, the idea now is to induct on ``levels" (each level is determined by the smallest index in the variables appearing in the matrix equations of the type above) supposing that all the variables appearing in the previous level (with the exception of variables of the form $d_{i}$ which must be treated separately later) have been demonstrated to be rational. At the $i$-th such level the arising matrix equation analogous to the one above is of the form:

\begin{eqnarray*}
{\scriptscriptstyle
\begin{psmatrix} 
1 & u_{k-i+1,k-i+2} & u_{k-i+1,k-i+3} & \hdots & u_{k-i+1,k} & v_{k-i+1,1} & v_{k-i+1,2} & \hdots & v_{k-i+1,m} \\
0 & 1 & u_{k-i+2,k-i+3} & \hdots & u_{k-i+2,k} & v_{k-i+2,1} & v_{k-i+2,2} & \hdots & v_{k-i+2,m}  \\
0 & 0 & 1 & \hdots & u_{k-i+3,k} & v_{k-i+3,1} & v_{k-i+3,2} & \hdots & v_{k-i+3,m}  \\
\vdots & \vdots & \vdots & \ddots & \vdots & \vdots & \vdots & \ddots & \vdots \\
0 & 0 & 0 & \hdots & 1 & v_{k,1} & v_{k,2} & \hdots & v_{k,m}  \\
p_{k-i+1,1} & p_{k-i+2,1} & p_{k-i+3,1} & \hdots & p_{k,1} & -1 & 0 & \hdots & 0 \\
p_{k-i+1,2} & p_{k-i+2,2} & p_{k-i+3,2} & \hdots & p_{k,2} & 0 & -1 & \hdots & 0 \\
\vdots & \vdots & \vdots & \ddots & \vdots & \vdots & \vdots & \ddots & \vdots \\
p_{k-i+1,m} & p_{k-i+2,m} & p_{k-i+3,m} & \hdots & p_{k,m} & 0 & 0 & \hdots & -1 \\
\end{psmatrix} \begin{psmatrix} u_{k-i,k-i+1} \\ u_{k-i,k-i+2} \\ u_{k-i,k-i+3} \\ \vdots \\ u_{k-i,k} \\ v_{k-i,1} \\ v_{k-i,2} \\ \vdots \\v_{k-i,m} \end{psmatrix} = \begin{psmatrix} q_{k-i+1,k-i} \\ q_{k-i+2,k-i} \\ q_{k-i+3,k-i} \\ \vdots \\ q_{k,k-i} \\ -p_{k-i,1} \\ -p_{k-i,2} \\ \vdots \\ -p_{k-i,m} \end{psmatrix}.}
\end{eqnarray*} 
\vspace{5 mm}

As all entries in the matrix on the left appear in the left or right hand side vector of some matrix equation from a previous level, the inductive hypothesis implies that they are rational. A few observations are in order. The first $i$ rows in the matrix above are linearly independent by the fact the first $i$ column sub-matrix is upper-triangular with ones along the diagonal. Second, the remaining $m$ row vectors have inner products with the first $i$ row vectors which are zero as the expressions resulting in computing these inner products come exactly as the equations $B_{0}\bp_{l}=\bwv_{l}$. Lastly, note that the last $m$ row vectors are linearly independent amongst themselves by the lower right block being minus the identity in $\real^{m\times m}$. Together, these observations justify the claim that the above matrix is invertible, so that the variables $u_{k-i,k-i+1},u_{k-i,k-i+2},\dots,u_{k-i,k},v_{k-i,1},\dots,v_{k-i,m}$ may be expressed as rationals. This completes the inductive portion of the argument. \\

Reflect on what is known about the variables which have appeared in this process so far. For each $i$, the variables $\{u_{k-i,j+1}\}_{j=k-i}^{k-1}$ have been shown to be rational along with the variables $\{v_{k-i,j}\}_{j=1}^m$. There is one set of equations which have not appeared yet, along with a set of variables which have yet to play a role (the variables $d_j$). Treating these will be the last step of this argument. \\

The diagonal elements of $BB^{\top}$ give rise to the equations
$$d_l^2 \left( 1+\sum\limits_{j=l}^{k-1} u_{l,j+1}^2+\sum\limits_{i=1}^{m} v_{l,i}^2 \right) = q_{l,l},\  l=1\dots,k,$$
where the convention is that a sum with starting index larger than the ending index is zero. For $l=k$, the corresponding equation is 
$$d_k^2 \left( 1+\sum\limits_{i=1}^{m} v_{k,i}^2 \right) = q_{k,k}.$$
Since all of the variables $v_{k,i}, q_{k,k}$ are rational, so is $d_k^2$. An analogous argument establishes that $d_l^2$ is rational as in those equations, $u_{l,j+1},q_{l,l}$ and $v_{l,j}$ are rational (by the previous inductive argument). All that remains now is to compute the inner products. These are of the form 
$$\langle \bwv_{i}, \bwv_{j} \rangle = \sum\limits_{l} d_{l}^2 v_{i,l}v_{j,l},\ \langle \bu_{i}, \bu_{j} \rangle=\sum\limits_{l} d_{l}^2 u_{i,l+1}u_{j,l+1},\ \langle \bu_{i}, \bwv_{j} \rangle=\sum\limits_{l}d_{l}^2 u_{i,l+1}v_{j,l},$$
which are all rational.
\endproof

We now show that conditions of Theorem~\ref{thm_rat} include tight frames, thus providing an alternate proof of Theorem~\ref{rational}.

\begin{cor} \label{cor_rat} Suppose that $B=\{ \bb_{i} \}_{i=1}^{n} \subset \real^{k}$ is a matrix with column vectors given by $\F$, a Parseval tight frame. Then $\spn_{\zed} \F$ is discrete if and only if $\langle \bb_{i}, \bb_{j} \rangle$ are rational.
\end{cor}

\proof
If the frame $\F$ is rational, then $\spn_{\zed} \F$ is a lattice by Proposition~1 of~\cite{etf2}. The reverse implication follows by setting $q_{i,j}=0$, $i\neq j$ and $q_{i,i}=1$ in Theorem~\ref{thm_rat} (after computing the QR decomposition of $B$).
\endproof

\proof[Second proof of Theorem~\ref{rational}] 
Suppose that $\F = \{ \bff_1,\dots,\bff_n \}$ is a uniform tight $(n,k)$-frame so that $L(\F) = \spn_{\zed} \F$ is a lattice. Then for all~$\bv \in \real^k$,
$$\|\bv\|^2 = \gamma \sum_{i=1}^n \langle \bv,\bff_i \rangle^2 = \sum_{i=1}^n \langle \bv, \sqrt{\gamma} \bff_i \rangle^2$$
for an appropriate constant~$\gamma > 0$. Hence $\F' = \sqrt{\gamma}\ \F$ is a Parseval tight frame and $\spn_{\zed} \F' = \sqrt{\gamma}\ L(\F)$ is again a lattice. Then Corollary~\ref{cor_rat} implies that inner products of vectors in $\F'$ are rational, and so inner products of vectors in $\F$ are rational multiples of~$1/\gamma$.
\endproof

\bigskip

\section{Lattices from irreducible group frames}
\label{grp_frames}

In this section we focus on group frames and lattices generated by them, in particular proving Theorem~\ref{Gframes}. As in Section~\ref{intro}, let $\bff_1 \in \real^k$ be a vector and let $G$ a finite group of orthogonal $k \times k$ matrices. Assume that 
$$\F := \left\{ U \bff_1 : U \in G \right\},$$
spans $\real^k$, that is, it is a $G$-frame. If $G$ is a cyclic group, $\F$ is called a {\it cyclic frame}. An example of a cyclic frame is the $(k,k+1)$-ETF discussed, for instance, in Section~5 of~\cite{etf1}:
\begin{equation}
\label{cyclic}
f_1=\frac{1}{\sqrt{k^2+k}} \left(\begin{array}{r} -k \\ 1 \\ \vdots \\ 1  \end{array}\right), \ldots, f_{k+1}=\frac{1}{\sqrt{k^2+k}} \left(\begin{array}{r} 1\\ 1\\ \vdots \\-k \end{array}\right).
\end{equation}
If $G$ is an abelian group, $\F$ is a {\it harmonic frame} (see Section~11.3 of~\cite{waldron_book}, Theorem~11.1). Notice that for any $G$-frame $\F$, $G \leq \Aut(\F)$. We also make a simple observation about the size of the $G$-frame~$\F$.

\begin{lem} \label{grp_size} Let $\F := \{ U \bff_1 : U \in G \}$ be a $G$-frame in~$\real^k$, then $|\F| = |G : G_{\bff_1}|$ where $G_{\bff_1}$ is the stabilizer of $\bff_1$ and $|\F| \leq |G|$. Further, $|\F| < |G|$ if and only if $\bff_1$ if an eigenvector for some non-identity matrix $W \in G$ with the corresponding eigenvalue equal to 1.
\end{lem}

\proof
The fact that $|\F| = |G : G_{\bff_1}| \leq |G|$ is clear from the definition. Now assume $|\F| < |G|$, which is equivalent to saying that $|G_{\bff_1}| > 1$. This is true if and only if there exists a non-identity matrix $W \in G$ such that $W\bff_1 = \bff_1$.
\endproof

\proof[Proof of Theorem~\ref{Gframes}]
The automorphism group of~$L(\F)$, $\Aut(L(\F))$, is the group of all orthogonal matrices that permute the lattice. Then we have
$$G \leq \Aut(\F) \leq \Aut(L(\F)),$$
and the action of $G$ on $\real^k$ is irreducible. Let $S(L(\F))$ be the set of minimal vectors of $L(\F)$ and let $E = \spn_\real S(L(\F))$. Since the automorphisms of $L(\F)$ permute the minimal vectors, it must be true that $E$ is closed under the action of $G$. Thus we must have $E = \real^k$, and so $G$ acts irreducibly on $E$, the space spanned by the minimal vectors of $L(\F)$. Then Theorem~3.6.6 of~\cite{martinet} guarantees that $S(L(\F))$ is a strongly eutactic configuration, and hence $L(\F)$ is a strongly eutactic lattice.
\endproof

\bigskip

\section{Vertex transitive graphs}
\label{graphs}

Construction of group frames from vertex transitive graphs is described in Section~10.7 of~\cite{waldron_book}\footnote{We found some comments to be misleading in this reference, such as in the proof of Proposition 10.2. That said, the treatment there is valuable, and overall the problems in this section are minor.}. We briefly review this subject here, proving Theorem~\ref{graph_thm} and providing some applications.

Let $\Gamma$ be a graph on $n$ vertices labeled by integers $1,\dots,n$ with automorphism group $G := \Aut(\Gamma)$. $\Gamma$ is called {\it vertex transitive} if for each pair of vertices $i,j$ there exists $\tau \in G$ such that $\tau(i)=j$. We define the {\it distance} between two vertices in a graph to be the number of edges in a shortest path connecting them. A connected graph $\Gamma$ is called {\it distance transitive} if for any two pairs of vertices $i,j$ and $k,l$ at the same distance from each other there existence an automorphism $\tau \in G$ such that $\tau(i) = k$ and $\tau(j) = l$. Clearly, distance transitive graphs are always vertex transitive, but the converse is not true. From here on graphs considered will always be vertex transitive, and we will indicate specifically when we need them to also be distance transitive. Let $\be_1,\dots,\be_n$ denote the standard basis vectors in $\real^n$. Then $G$ acts on $\real^n$ by
$$\tau \left( \sum_{i=1}^n c_i \be_i \right) = \sum_{i=1}^n c_i \be_{\tau(i)}$$
for every $\tau \in G$ and vector $\sum_{i=1}^n c_i \be_i \in \real^n$. Let $A = (a_{ij})$ be the $n \times n$ adjacency matrix of $\Gamma$, so that $a_{ij} = 1$ if vertices $i$ and $j$ are connected by an edge and $a_{ij}=0$ otherwise. Then $a_{\tau(i) \tau(j)}=a_{ij}$ for all $\tau\in G$. The matrix $A$ is symmetric, with real eigenvalues $\lambda_1,\dots,\lambda_k$, each of multiplicity $m_{\lambda_i}$, so that $\sum_{i=1}^k m_{\lambda_i} = n$. From now on, we call these the eigenvalues of the graph~$\Gamma$. For each $\lambda_i$ let $V_{\lambda_i} \subset \real^n$ be the corresponding $m_{\lambda_i}$-dimensional eigenspace. The group $G$ acts on each eigenspace $V_{\lambda_i}$ and for any nonzero vector $\bv \in V_{\lambda_i}$ the orbit $G\bv$ of $\bv$ under the action of $G$ is a group frame in $V_{\lambda_i} \cong \real^{m_i}$. When $\Gamma$ is a distance transitive graph, this action of $G$ on $V_{\lambda_i}$ is irreducible, hence producing an irreducible group frame (see Proposition~4.1.11 on p.~137 of~\cite{distance_regular}). Further, if $P_{\lambda_i}$ is the orthogonal projection onto $V_{\lambda_i}$, then for any $\tau \in G$ and $\bx \in \real^n$,
$$\tau (P_{\lambda_i}(\bx)) = P_{\lambda_i} (\tau(\bx)).$$
As indicated in Section~10.7 of~\cite{waldron_book}, this identity is true since the action of $\tau \in G$ and the action of the adjacency matrix $A$ on a vector commute, i.e.
$$\tau(A \be_k)=\sum \limits_i a_{ik} \tau \be_i = \sum \limits_j a_{\tau^{-1}(j),k} \be_j=\sum\limits_j a_{j, \tau(k)} \be_j=A(\tau(\be_k)).$$

\proof[Proof of Theorem~\ref{graph_thm}]
Suppose now that an eigenvalue $\lambda_i$ is an integer. We know that the group $G$ consists of permutation matrices. Pick a nonzero integer vector $\bx \in \real^n$. Then $P_{\lambda_i} \bx \in V_{\lambda_i}$ and the frame $\F_{\lambda_i}(\bx) := G P_{\lambda_i} \bx = P_{\lambda_i} (G \bx)$ is rational, and hence generates a lattice $L(\F_{\lambda_i}(\bx)) = \spn_{\zed} \F_{\lambda_i}(\bx)$. This lattice is strongly eutactic if this group frame is irreducible, which is the case when the graph is distance transitive. Let $H$ be the kernel of the action of $G$ on $V_{\lambda_i}$, i.e. 
$$H = \left\{ \tau \in G : \tau(\bx) = \bx\ \text{ for all }\ \bx \in V_{\lambda_i} \right\}.$$
Notice that $H$ is a normal subgroup of $G$, since for any $\sigma \in G$ and $\bx \in V_{\lambda_i}$, $\sigma(\bx) \in V_{\lambda_i}$, and so
$$(\tau \sigma) (\bx) = \tau (\sigma(\bx)) = \sigma(\bx) = \sigma (\tau(\bx)) = (\sigma \tau) (\bx),$$
for any $\tau \in H$. Then the quotient group $G/H$ is isomorphic to a subgroup of $\Aut(L(\F_{\lambda_i}))$. If $\bx = \be_1$, then the corresponding frame
$$\F_{\lambda_i} := \F_{\lambda_i}(\be_1)$$
consists of column vectors of $P_{\lambda_i}$ (possibly with repetitions), since $\tau \be_1$ is some $\be_j$ for every $\tau \in G$, and every $\be_j$ is representable as $\tau \be_1$ for some $\tau \in G$, since the graph is vertex transitive. Then the resulting lattice $L(\F_{\lambda_i}) = P_{\lambda_i} \zed^n$, and this concludes the proof of Theorem~\ref{graph_thm}. 
\endproof

\noindent
We refer to the lattice $L(\F_{\lambda_i})$ described above as lattice {\it generated} by the graph $\Gamma$ and denote it by $L_{\Gamma, \lambda_i}$.

\begin{rem} \label{vertex_vs_distance} While our proof that the lattice $L_{\Gamma, \lambda_i}$ is strongly eutactic only applies to the situations when $\Gamma$ is distance transitive, there are examples of vertex transitive graphs which are not distance transitive that nonetheless still produce strongly eutactic lattices: we demonstrate some such examples below. It would be interesting to understand if this is indeed the case for all vertex transitive graphs, or if there exist some that generate lattices that are not strongly eutactic.
\end{rem}

For the rest of this section, we consider examples of this lattice construction when applied to various graphs and their products. One class of lattices that will figure prominently in our examples are {\it root lattices}, that is, integral lattices generated by vectors of norm~$2$, which are called its {\it roots} (recall that a lattice is integral if the inner product between any two vectors is always an integer). Also recall that the {\it dual lattice} of a full rank lattice $L \subset \real^n$ is
$$L^* := \left\{ \bx \in \real^n : \left< \bx, \bwy \right> \in \zed \text{ for all } \bwy \in L \right\}.$$
If $L$ is integral, then $L \subseteq L^*$.

\begin{lem} \label{disconnect} Let $0_n$ be a completely disconnected graph on $n$ vertices, then $0_n$ generates the integer lattice~$\zed^n$.
\end{lem}

\proof
The adjacency matrix for $0_n$ is the $n \times n$ $0$-matrix, and so it has one eigenvalue $0$ with multiplicity $n$ with the corresponding eigenspace being the entire~$\real^n$. The automorphism group of $0_n$ is $S_n$, so the group frame obtained from the vector~$\be_1$ is the full standard basis, which spans the lattice~$\zed^n$.
\endproof

\begin{lem} \label{complete} The complete graph $K_n$ generates (a lattice similar to) the root lattice
$$A_{n-1} = \left\{ \bx \in \zed^n : \sum_{i=1}^n x_i = 0 \right\}.$$
\end{lem}

\proof
The complete graph $K_n$ is the graph on $n$ vertices with no loops in which every vertex is connected to every other. Hence adjacency matrix $A$ has $1$'s for all the off-diagonal entries and $0$'s on the diagonal. There are two eigenvalues: $\lambda_1=-1$ with multiplicity $n-1$ and $\lambda_2=n-1$ with multiplicity $1$. The eigenspace corresponding to $\lambda_2$ is $V_{n-1} = \spn_{\real} \{ (1,\dots,1)^\top \}$ and the eigenspace $V_{-1}$ corresponding to $\lambda_1$ is the orthogonal complement of $V_{n-1}$ in~$\real^n$. The automorphism group of $K_n$ is $S_n$. The orthogonal projection onto $V_{-1}$ is given by
$$P_{-1} = \frac{1}{n-1} \begin{pmatrix} n-1 & -1 & \dots & -1 \\ -1 & n-1 & \dots & -1 \\ \vdots & \vdots & \ddots & \vdots \\ -1 & -1 & \dots & n-1 \end{pmatrix},$$
so the lattice $L_{K_n,-1}$ generated by the columns of $P_{-1}$ is the root lattice $A_{n-1} = \zed^n \cap V_{-1}$ rescaled by the factor~$1/(n-1)$.
\endproof

Next we consider graphs that are constructed as products of smaller graphs. We start with disjoint unions. In order for such a graph to be vertex transitive, all the components in the disjoint union need to be vertex transitive and isomorphic to each other. Hence we can think of them as copies of the same vertex transitive graph.

\begin{lem} \label{disjoint} Let $\Gamma$ be a vertex transitive graph constructed as a disjoint union of $k$ copies of a vertex transitive graph $\Delta$. Let $\lambda$ be a rational eigenvalue of $\Delta$ and $L_{\Delta,\lambda}$ be a lattice generated by the $\lambda$-eigenspace of $\Delta$. Then $\Gamma$  also has $\lambda$ as an eigenvalue and generates a lattice given by the orthogonal sum of $k$ copies of~$L_{\Delta,\lambda}$.
\end{lem}

\proof
Let $m$ be the number of vertices of $\Delta$ and let $A_{\Delta}$ be its adjacency matrix. Then the $mk \times mk$ adjacency matrix $A_{\Gamma}$ of the graph $\Gamma$ is a block matrix with diagonal $m \times m$ blocks being $A_{\Delta}$ and the rest filled up with $0$ blocks, i.e.
$$A_{\Gamma} = \begin{pmatrix} A_{\Delta} & 0 & \dots & 0 \\ 0 & A_{\Delta} & \dots & 0 \\ \vdots & \vdots & \ddots & \vdots \\ 0 & 0 & \dots & A_{\Delta} \end{pmatrix}.$$
Let us refer to a block matrix like this as $\bigoplus_k(A_{\Delta})$. $A_{\Gamma}$ has the same eigenvalues as $A_{\Delta}$, but of $k$ times greater multiplicity. Let $V_{\Delta,\lambda}$ be the $\lambda$-eigenspace of $A_{\Delta}$ with the corresponding projection matrix $P_{\Delta,\lambda}$. The $\lambda$-eigenspace of $A_{\Gamma}$ is the orthogonal sum of $k$ copies of $V_{\Delta,\lambda}$ and the corresponding projection matrix is $\bigoplus_k(P_{\Delta,\lambda})$. Hence the lattice $L_{\Gamma,\lambda}$ generated by the column vectors of this matrix is the orthogonal sum of $k$ copies of~$L_{\Delta,\lambda}$.
\endproof
\smallskip

\noindent
Now we recall the three fundamental commutative product constructions of graphs (see~\cite{products} and~\cite{cartesian} for detailed information). In each of these constructions, each eigenvalue $\nu$ of the product graph $\Gamma$ is derived from a pair of eigenvalues $\lambda$ and $\mu$ of the component graphs $\Delta_1$ and $\Delta_2$, respectively, via some function $f(\lambda,\mu)$. This function $f$ differs depending on which product we consider. Spectral properties of product graphs are nicely summarized in~\cite{sayama}.
\smallskip

The {\it Cartesian product} of two graphs $\Delta_1$ and $\Delta_2$, denoted $\Delta_1 \square \Delta_2$, is the graph whose vertices are pairs $(u,v)$, where $u$ is a vertex of $\Delta_1$ and $v$ is a vertex of $\Delta_2$, and two vertices $(u_1,v_1)$ and $(u_2,v_2)$ are connected by an edge if and only if either $u_1=u_2$ and $v_1,v_2$ are connected by an edge in $\Delta_2$, or $v_1=v_2$ and $u_1,u_2$ are connected by an edge in $\Delta_1$. Then  $\Delta_1 \square \Delta_2$ is vertex transitive if and only if both $\Delta_1$ and $\Delta_2$ are vertex transitive (\cite{godsil}, Section~7.14, or~\cite{products}). For each pair of eigenvalues $\lambda$ of $\Delta_1$ and $\mu$ of $\Delta_2$, there is an eigenvalue $\nu$ of $\Delta_1 \square \Delta_2$ given by 
$$\nu = f(\lambda,\mu) := \lambda + \mu,$$
and if $\bu, \bv$ are corresponding eigenvectors of $\Delta_1, \Delta_2$, respectively, then $\bu \otimes \bv$ is an eigenvector of $\Gamma$ corresponding to~$\nu$.
\smallskip

The {\it direct product} of two graphs $\Delta_1$ and $\Delta_2$, denoted $\Delta_1 \times \Delta_2$ is the graph whose vertices are pairs $(u,v)$, where $u$ is a vertex of $\Delta_1$ and $v$ is a vertex of $\Delta_2$, and two vertices $(u_1,v_1)$ and $(u_2,v_2)$ are connected by an edge if and only if both pairs $u_1, u_2$ and $v_1, v_2$ are connected by an edge in $\Delta_1, \Delta_2$, respectively. If $\Delta_1$ and $\Delta_2$ are vertex transitive, then $\Delta_1 \times \Delta_2$ is vertex transitive. The converse statement is not as straight-forward, and distinguishes between bipartite and non-bipartite graphs (see \cite{direct_products}). For each pair of eigenvalues $\lambda$ of $\Delta_1$ and $\mu$ of $\Delta_2$, there is an eigenvalue $\nu$ of $\Delta_1 \times \Delta_2$ given by 
$$\nu = f(\lambda,\mu) := \lambda \mu,$$
and if $\bu, \bv$ are corresponding eigenvectors of $\Delta_1, \Delta_2$, respectively, then $\bu \otimes \bv$ is an eigenvector of $\Gamma$ corresponding to~$\nu$.
\smallskip

The {\it strong product} of two graphs $\Delta_1$ and $\Delta_2$, denoted $\Delta_1 \boxtimes \Delta_2$, is the graph whose vertices are pairs $(u,v)$, where $u$ is a vertex of $\Delta_1$ and $v$ is a vertex of $\Delta_2$, and two vertices $(u_1,v_1)$ and $(u_2,v_2)$ are connected by an edge if and only if $u_1, u_2$ and $v_1, v_2$ are either equal or connected by an edge in $\Delta_1, \Delta_2$, respectively. The graph  $\Delta_1 \boxtimes \Delta_2$ is vertex transitive if and only if both $\Delta_1$ and $\Delta_2$ are vertex transitive (Section~7.4 of~\cite{products}). For each pair of eigenvalues $\lambda$ of $\Delta_1$ and $\mu$ of $\Delta_2$, there is an eigenvalue $\nu$ of $\Delta_1 \boxtimes \Delta_2$ given by 
$$\nu = f(\lambda,\mu) := (\lambda + 1)(\mu + 1) - 1,$$
and if $\bu, \bv$ are corresponding eigenvectors of $\Delta_1, \Delta_2$, respectively, then $\bu \otimes \bv$ is an eigenvector of $\Gamma$ corresponding to~$\nu$.
\smallskip

The {\it lexicographic product} of two vertex transitive graphs $\Delta_1$ and $\Delta_2$ is a vertex transitive graph whose vertices are pairs $(u,v)$, where $u$ is a vertex of $\Delta_1$ and $v$ is a vertex of $\Delta_2$, and two vertices $(u_1,v_1)$ and $(u_2,v_2)$ are connected by an edge if and only if either $u_1$, $u_2$ are connected in $\Delta_1$, or $u_1=u_2$ and $v_1,v_2$ are connected in $\Delta_2$.
\smallskip

For two vectors $\bx \in \real^{m_1}, \bwy \in \real^{m_2}$ and $m_1 \times m_1$, $m_2 \times m_2$ matrices $A, B$, respectively, we have
\begin{equation}
\label{tensor_eq}
(A \bx) \otimes (B \bwy) = (A \otimes B) (\bx \otimes \bwy),
\end{equation}
where $\otimes$ stands for the usual Kronecker (outer) product of matrices and vectors. Further, if two vectors $\bx_1,\bx_2 \in \real^{m_1}$ are orthogonal and $\bwy \in \real^{m_2}$, then simple tensors $\bx_1 \otimes \bwy$ and $\bx_2 \otimes \bwy$ are also orthogonal.

\begin{thm} \label{tensor} Let $\Delta_1,\Delta_2$ be vertex transitive graphs on $m_1$, $m_2$ vertices, respectively, and let $\Gamma$ be a product graph
$$\Gamma = \Delta_1 * \Delta_2$$
on $m_1 m_2$ vertices, where $*$ stands for $\square$, $\times$, or $\boxtimes$. Let $\nu$ be an eigenvalue of $\Gamma$ and $(\lambda_i,\mu_i)$ for $1 \leq i \leq k$ pairs of eigenvalues of $\Delta_1,\Delta_2$ respectively so that
$$\nu = f(\lambda_i,\mu_i)\ \text{for all }\ 1 \leq i \leq k$$
for the appropriate $f$. Let $L_{\Delta_1,\lambda_i}$ and $L_{\Delta_2,\mu_i}$ for each $1 \leq i \leq k$ be the corresponding lattices. Then $L_{\Gamma,\nu} $ is the orthogonal projection of $\zed^{m_1 m_2}$ onto the space spanned by
$$\left( L_{\Delta_1,\lambda_1} \otimes_{\zed} L_{\Delta_2,\mu_1} \right) \oplus \cdots \oplus  \left( L_{\Delta_1,\lambda_k} \otimes_{\zed} L_{\Delta_2,\mu_k} \right),$$
where $\oplus$ is the orthogonal direct sum. In particular, if $k=1$ then
$$L_{\Gamma,\nu} = L_{\Delta_1,\lambda_1} \otimes_{\zed} L_{\Delta_2,\mu_1},$$
up to similarity.
\end{thm} 

\proof
Let $V_{\Delta_1,\lambda_i}$, $W_{\Delta_2,\mu_i}$ be the eigenspaces of $\Delta_1$, $\Delta_2$ corresponding to $\lambda_i$, $\mu_i$, respectively, with the corresponding orthogonal projection matrices $P_{\Delta_1,\lambda_i}$, $P_{\Delta_2,\mu_i}$. Then
$$L_{\Delta_1,\lambda_i} = P_{\Delta_1,\lambda_i} \zed^{m_1} \subset V_{\Delta_1,\lambda_i},\ L_{\Delta_2,\mu_i} = P_{\Delta_2,\mu_i} \zed^{m_2} \subset W_{\Delta_2,\mu_i},$$
and $V_{\Delta_1,\lambda_i} = \spn_{\real} L_{\Delta_1,\lambda_i}$, $W_{\Delta_2,\mu_i} = \spn_{\real} L_{\Delta_2,\mu_i}$, so
$$V_{\Delta_1,\lambda_i} \otimes_{\real} W_{\Delta_2,\mu_i} = \spn_{\real} \left( L_{\Delta_1,\lambda_i} \otimes_{\zed} L_{\Delta_2,\mu_i} \right).$$
Since adjacency matrices of graphs are symmetric, the eigenspaces corresponding to distinct eigenvalues are orthogonal, so that any two $V_{\Delta_1,\lambda_i}$ are orthogonal to each other, as are any two~$W_{\Delta_2,\mu_i}$. Then each two $V_{\Delta_1,\lambda_i} \otimes_{\real} W_{\Delta_2,\mu_i}$ are also orthogonal to each other, and the eigenspace of~$\Gamma$ corresponding to~$\nu$ is
\begin{eqnarray*}
U_{\Gamma,\nu} = P_{\Gamma,\nu} \real^{m_1 m_2} & = & \left( P_{\Delta_1,\lambda_1} \otimes P_{\Delta_2,\mu_1} \right) \real^{m_1 m_2} \oplus \cdots \oplus \left( P_{\Delta_1,\lambda_k} \otimes P_{\Delta_2,\mu_k} \right) \real^{m_1 m_2} \\
& = & \left( P_{\Delta_1,\lambda_1} \real^{m_1}  \otimes_{\real} P_{\Delta_2,\mu_1} \real^{m_2} \right) \oplus \cdots \oplus \left( P_{\Delta_1,\lambda_k} \real^{m_1}  \otimes_{\real} P_{\Delta_2,\mu_k} \real^{m_2} \right) \\
& = & \left( V_{\Delta_1,\lambda_1} \otimes_{\real} W_{\Delta_2,\mu_1} \right) \oplus \cdots \oplus  \left( V_{\Delta_1,\lambda_k} \otimes_{\real} W_{\Delta_2,\mu_k} \right),
\end{eqnarray*}
by~\eqref{tensor_eq}, where $P_{\Gamma,\nu}$ is the orthogonal projection matrix onto $U_{\Gamma,\nu}$; we are using here the fact that $\real^{m_1} \otimes_{\real} \real^{m_2} = \real^{m_1 m_2}$. Then $L_{\Gamma,\nu} = P_{\Gamma,\nu} \zed^{m_1 m_2}$.

Now suppose $k=1$, then applying~\eqref{tensor_eq} again and using the fact that $\zed^{m_1} \otimes_{\zed} \zed^{m_2} = \zed^{m_1 m_2}$, we have:
$$L_{\Gamma,\nu} = P_{\Gamma,\nu} \zed^{m_1 m_2} = \left( P_{\Delta_1,\lambda_1} \otimes P_{\Delta_2,\mu_1} \right) \zed^{m_1 m_2} = P_{\Delta_1,\lambda_1} \zed^{m_1} \otimes_{\zed} P_{\Delta_2,\mu_1} \zed^{m_2} = L_{\Delta_1,\lambda_1} \otimes_{\zed} L_{\Delta_2,\mu_1}.$$
This completes the proof.
\endproof

\noindent
\begin{ex} \label{product_ex} Let $\Delta_1$ be the complete graph $K_3$ and $\Delta_2$ the $4$-cycle graph~$C_4$. Eigenvalues of $K_3$ are $\lambda_1 = 2$ (multiplicity $1$) and $\lambda_2 = -1$ (multiplicity $2$); eigenvalues of $C_4$ are $\mu_1 = 2$ (multiplicity $1$), $\mu_2 = -2$ (multiplicity $1$), $\mu_3 = 0$ (multiplicity $2$). The corresponding lattices are
$$L_{K_3,2} = \frac{1}{3} \begin{pmatrix} \ 1 \\ \ 1 \\ \ 1 \end{pmatrix} \zed,\ L_{K_3,-1} = \frac{1}{3} \begin{pmatrix} \ \ 2 & -1 \\ -1 & \ \ 2 \\ -1 & -1 \end{pmatrix} \zed^2,$$
and
$$L_{C_4,2} = \frac{1}{4} \begin{pmatrix} \ 1 \\ \ 1 \\ \ 1 \\ \ 1 \end{pmatrix} \zed,\ L_{C_4,-2} = \frac{1}{4} \begin{pmatrix} \ \ 1 \\ -1 \\ \ \ 1 \\ -1 \end{pmatrix} \zed,\ L_{C_4,0} = \frac{1}{4} \begin{pmatrix} \ 1 & \ \ 0 \\ \ 0 &\ \ 1 \\ -1 & \ \ 0 \\ \ \ 0 & -1 \end{pmatrix} \zed^2.$$
\smallskip

\noindent
Let $\Gamma_1 = K_3 \square C_4$, then $\nu = -1$ is an eigenvalue of $\Gamma_1$, obtained in a unique way as $\nu = \lambda_2 + \mu_3$,
hence
$$L_{\Gamma_1,-1} = L_{K_3,-1} \otimes_{\zed} L_{C_4,0} \sim A_2 \otimes _{\zed} \zed^2 = A_2 \oplus A_2.$$

\noindent
Let $\Gamma_2 = K_3 \times C_4$, then $\nu = 0$ is an eigenvalue of $\Gamma_2$, obtained as 
$$\nu = \lambda_1 \mu_3 = \lambda_2 \mu_3,$$
hence $L_{\Gamma_2,0}$ is the orthogonal projection of $\zed^{12}$ onto the space spanned by
$$\left( L_{K_3,2} \otimes_{\zed} L_{C_4,0} \right) \oplus \left( L_{K_3,-1} \otimes_{\zed} L_{C_4,0} \right) = \left( L_{K_3,2} \oplus L_{K_3,-1} \right) \otimes_{\zed} L_{C_4,0} \sim \zed^3 \otimes_{\zed} \zed^2 = \zed^6.$$
Hence $L_{\Gamma_2,0}$ is similar to~$\zed^6$.
\smallskip

\noindent
Let $\Gamma_3 = K_3 \boxtimes C_4$, then $\nu = -1$ is an eigenvalue of $\Gamma_2$, obtained as 
$$\nu = (\lambda_1 + 1)( \mu_1 + 1) - 1 =(\lambda_1 + 1)( \mu_2 + 1) - 1 = (\lambda_1 + 1)( \mu_3 + 1) - 1,$$
hence $L_{\Gamma_3,-1}$ is the orthogonal projection of $\zed^{12}$ onto the space spanned by
\begin{eqnarray*}
& & \left( L_{K_3,-1} \otimes_{\zed} L_{C_4,2} \right) \oplus \left( L_{K_3,-1} \otimes_{\zed} L_{C_4,-2} \right) \oplus \left( L_{K_3,-1} \otimes_{\zed} L_{C_4,0} \right)  \\
& = &  L_{K_3,-1} \otimes_{\zed} \left( L_{C_4,2} \oplus L_{C_4,-2} \oplus L_{C_4,0} \right) \sim A_2 \otimes_{\zed} \zed^4 \\
& = & A_2 \oplus A_2 \oplus A_2 \oplus A_2.
\end{eqnarray*}
Hence $L_{\Gamma_2,0}$ is similar to~$A_2 \oplus A_2 \oplus A_2 \oplus A_2$.
\smallskip

\noindent
Let $\Gamma_4 = K_3 \circ C_4$ be the lexicographic product of $K_3$ by $C_4$. Unlike the previously considered products, this one is not commutative.Then eigenvalues of $\Gamma_4$ are $10$ (multiplicity $1$), $0$ (multiplicity $6$), $-2$ (multiplicity $5$). The lattice $L_{\Gamma_4,-2}$ is similar to~$A_5^*$, and the lattice $L_{\Gamma_4,0}$ is similar to~$\zed^6$.
\end{ex}
\medskip

We also discuss a relation between lattices generated by a graph and by its complement. If $\Gamma$ is a graph on $n$ vertices, then its complement $\Gamma'$ is a graph on the same vertices that has no common edges with $\Gamma$ and so when `put together' the two form a complete graph $K_n$. Vertex transitive graphs are regular, so let $k$ be the common degree of the vertices of $\Gamma$. Then $n-k-1$ is the common degree of the vertices of $\Gamma'$. So $k$ is an eigenvalue of $\Gamma$ of multiplicity $1$ with the corresponding eigenvector ${\bf 1}:=(1,\dots,1)^{\top}$ and $n-k-1$ is an eigenvalue of $\Gamma'$ of the same multiplicity with the same corresponding eigenvector. Moreover the following result holds.

\begin{prop} \label{complement} Let $\Gamma$ be a vertex transitive graph on $n$ vertices of degree $k$ and $\Gamma'$ its complement. Then for each eigenvalue $\lambda \neq k$ of $\Gamma$ there is an eigenvalue $\lambda' = -\lambda-1$ of $\Gamma'$ of the same multiplicity and the lattices $L_{\Gamma,\lambda}$ and $L_{\Gamma',\lambda'}$ are the same.
\end{prop}

\proof
It is well known that if $p(x)$ is the characteristic polynomial of the adjacency matrix $A$ of $\Gamma$, then the characteristic polynomial of the adjacency matrix $B$ of $\Gamma'$ is
$$q(x) = (-1)^n \frac{x-n+k+1}{x+k+1} p(-x-1),$$
and so for each eigenvalue $\lambda \neq k$ of $\Gamma$ there is an eigenvalue $\lambda' = -\lambda-1$ of $\Gamma'$ of the same multiplicity (see, for instance, p.~27 of~\cite{spectra}). Further, the adjacency matrices satisfy the relation
$$B = J_n - I_n - A,$$
where $I_n$ is the $n \times n$ identity matrix and $J_n$ is the $n \times n$ matrix consisting of all~$1$'s. Let $\lambda \neq k$ be an eigenvalue of $\Gamma$ with a corresponding eigenvector $\bx$. Since eigenspaces of $\Gamma$ corresponding to different eigenvalues are orthogonal, $\bx$ must be orthogonal to ${\bf 1}$, which means that
$$\sum_{i=1}^n x_i = 0,$$
and so $J_n \bx = \bo$. Then
$$B \bx = J_n \bx - \bx - \lambda \bx = (-\lambda - 1) \bx,$$ 
i.e. $\bx$ is an eigenvector of $B$ corresponding to the eigenvalue $\lambda'$. This means that the eigenspace of $\Gamma'$ corresponding to the eigenvalue $\lambda' = - \lambda -1$ is the same as the eigenspace of $\Gamma$ corresponding to the eigenvalue $\lambda$, hence they generate the same lattices.
\endproof

We now consider more examples. In all the examples to follow, lattices are specified up to similarity. Information about the graphs we mention can be found, for instance, in~\cite{spectra}.

\begin{ex} \label{hamming} Recall the construction of the Hamming graph $H(d,q)$: if $S$ is a set of $q$ elements and $d$ a positive integer, then vertex set of $H(d,q)$ is $S^d$, the set of ordered $d$-tuples of elements of $S$, and two vertices are connected by an edge if they differ in precisely one coordinate. $H(d,q)$ has eigenvalues $(q-1)d-qi$ with multiplicity $\binom{d}{i} (q-1)^i$ for $0 \leq i \leq d$. It is well known that $H(d,q)$ is the Cartesian product of $d$ complete graphs $K_q$, and hence gives rise to product lattices. Hamming graphs are known to be distance transitive.

For instance, $H(2,3)$ has 9 vertices and three eigenvalues: $4$ (multiplicity $1$), $-2$ (multiplicity $4$) and $1$ (multiplicity $4$). Projection matrices of both of the $4$-dimensional eigenspaces give rise to the same tensor product lattice: $A_2 \otimes_{\zed} A_2$.

On the other hand, the graph $H(3,2)$ has $8$ vertices and is isomorphic to the cube graph~$Q_3$, i.e. 
$$H(3,2) = K_2 \square K_2 \square K_2 = K_2 \square C_4,$$
where $C_4$ is as in Example~\ref{product_ex} with eigenvalues $\mu_1$, $\mu_2$, $\mu_3$ and the corresponding lattices, and $K_2$ that has multiplicity $1$ eigenvalues $\lambda_1 = 1$, $\lambda_2 = -1$ with
$$L_{K_2,1} = \frac{1}{2} \begin{pmatrix} 1 \\ 1 \end{pmatrix} \zed,\ L_{K_2,-1} = \frac{1}{2} \begin{pmatrix} 1 \\ -1 \end{pmatrix} \zed.$$
Therefore eigenvalues of $H(3,2)$ are:
\begin{itemize}
\item $1$ (multiplicity $3$), obtained in $2$ ways: $\lambda_1 + \mu_3 = 1 + 0$ and $\lambda_2 + \mu_1 = -1 + 2$;
\item $-1$ (multiplicity $3$), obtained in $2$ ways: $\lambda_1 + \mu_2 = 1 + (-2)$ and $\lambda_2 + \mu_3 = -1 + 0$;
\item $3$ (multiplicity $1$), obtained as $\lambda_1+\mu_1$;
\item $-3$ (multiplicity $1$), obtained as $\lambda_2+\mu_2$.
\end{itemize}
The lattices $L_{H(3,2),3}$ and $L_{H(3,2),3}$ are both similar to~$\zed$, however $L_{H(3,2),1}$ is the orthogonal projection of $\zed^{8}$ onto the space spanned by
$$\left( L_{K_2,1} \otimes_{\zed} L_{C_4,0} \right) \oplus \left( L_{K_2,-1} \otimes_{\zed} L_{C_4,2} \right).$$
This lattice is similar to~$A_3^*$, and the same is true for the lattice~$L_{H(3,2),-1}$. This example demonstrates that a product graph construction can generate a lattice that is not a tensor product or direct sum.
\end{ex}

\begin{ex} \label{petersen} Recall the construction of the Kneser graph $KG_{n,k}$: vertices of this graph correspond to $k$-element subsets of a set of $n$ elements, and two vertices are connected by an edge if the corresponding sets are disjoint. $KG_{n,k}$ has eigenvalue $(-1)^j \binom{n-k-j}{k-j}$ occurring with multiplicity $\binom{n}{j} - \binom{n}{j-1}$ for all $j = 1,\dots,k$, and therefore gives rise to lattices in arbitrarily large dimensions. While Kneser graphs are not distance transitive in general, there are some examples that are.

For instance, Petersen graph (which is the same as the Kneser graph $KG_{5,2}$) has $10$ vertices and three eigenvalues: $3$ (multiplicity 1), $1$ (multiplicity $5$) and $-2$ (multiplicity) $4$. It is distance transitive, and hence generates strongly eutactic lattices corresponding to its eigenvalues. For eigenvalue $-2$, we obtain the lattice $A_4^*$. For eigenvalue $1$, we obtain $A_5^2$, an example of the Coxeter-Barnes lattice $A_n^r$, defined as the lattice contained in the hyperplane $H=(\be_{1}+\dots+\be_{n+1})^{\perp}$ with the basis 
$$\left\{ \be_{1}-\be_{2},\dots,\be_{1}-\be_{n},\frac{1}{r} \sum\limits_{i=2}^{n+1} (\be_{1}-\be_{i}) \right\}$$
and defined for all positive rational $r$. When $r$ is an integer dividing $n+1$, these are exactly the lattices $\Lambda$ for which $A_n \subset \Lambda \subset A_n^{*}$, so that $A_n^{r}$ contains $A_n$ to index $r$ (\cite{martinet}, Section~5.2). In particular, $A_5^2$ is the unique sublattice of the dual lattice
$$A_5^* := \{ \bx \in \real^5 : \bx^{\top} \bwy \in \zed\ \text{for all}\ \bwy \in A_5 \},$$
which contains $A_5$ to index 2. As mentioned above, it can be described as a full rank lattice in the hyperplane
$$\left\{ \bx \in \real^6 : \sum_{i=1}^6 x_i = 0 \right\},$$
identified with~$\real^5$. Here is this description:
$$A_5^2 = \spn_{\zed} \left\{ \be_1 - \be_2, \dots, \be_1 - \be_5, \frac{1}{2} \left( 5\be_1 - \sum_{i=2}^6 \be_i \right) \right\},$$
where $\be_1,\dots,\be_6$ are standard basis vectors in~$\real^6$. 
\end{ex}

\begin{ex} \label{line} The line graph of a graph $\Gamma$ is the graph $\Gamma'$ whose vertices correspond to edges of $\Gamma$, and two vertices are connected by an edge if and only if the corresponding edges in $\Gamma$ meet in a vertex. For instance, the line graph of the Petersen graph is a distance transitive graph on $15$ vertices. Among its eigenvalues, $-1$ comes with multiplicity $4$ and the corresponding lattice is $A_4^*$, $-2$ comes with multiplicity $5$ and the corresponding lattice is the Coxeter lattice $A_5^3$, which can be described as a full rank lattice in the hyperplane
$$\left\{ \bx \in \real^6 : \sum_{i=1}^6 x_i = 0 \right\},$$
identified with~$\real^5$. Here is the description:
$$A_5^3 = \spn_{\zed} \left\{ \be_1 - \be_2, \dots, \be_1 - \be_5, \frac{1}{3} \left( 5\be_1 - \sum_{i=2}^6 \be_i \right) \right\},$$
where $\be_1,\dots,\be_6$ are standard basis vectors in~$\real^6$. It is the unique sublattice of $A_5^*$ containing $A_5$ to index~$3$; it is isometric to the dual of~$A_5^2$.
\end{ex}

\begin{ex} \label{johnson} Recall the construction of the Johnson graph $J(n,k)$: vertices of this graph correspond to $k$-element subsets of a set of $n$ elements, and two vertices are connected by an edge if the corresponding sets intersect in $k-1$ elements. $J(n,k)$ is a distance transitive graph, which has $\binom{n}{k}$ vertices and  eigenvalue $((k-j)(n-k-j)-j)$ occurring with multiplicity $\binom{n}{j} - \binom{n}{j-1}$ for all $j = 1,\dots,\min \{ k,n-k \}$, and therefore gives rise to strongly eutactic lattices in arbitrarily large dimensions.

It is well known that Johnson graph $J(n,2)$ (also known as the triangular graph~$T_n$) is the line graph of the complete graph $K_n$ and the complement of the Kneser graph $KG_{n,2}$. In particular, $J(5,2)$ is the line graph of $K_5$ and the complement of the Petersen graph. Further, $J(n,2)$ is a strongly regular graph, and so always has three eigenvalues: $2(n-2)$ (multiplicity $1$), $n-4$ (multiplicity $n-1$), $-2$ (multiplicity $n(n-3)/2$). We present some examples of lattices from $J(n,2)$ in Table~\ref{J_graph_ex}, which are the same as for its complement $KG_{n,2}$. In this table, the lattice $L_{J(n,2),-2}$ for $n=6$ is listed as the $9$-dimensional lattice {\tt sth15} in the online catalog~\cite{streut} of strongly eutactic lattices; for larger $n$ in our table these lattices are not catalogued.

\begin{center} 
\begin{table}[!ht]
\begin{tabular}{|P{2.5cm}|P{2.5cm}|P{1.5cm}|P{3.5cm}|} \hline
{\bf $J(n,2)$} & {\bf \# of vertices} & {\bf Lattice $L_{J(n,2),n-4}$} & {\bf Lattice $L_{J(n,2),-2}$} \\
\hline \hline
$J(4,2)$ & $(6)$ & $\zed^3$ & $A_2$ \\ \hline
$J(5,2)$ & $(10)$ & $A_4^*$ & $A_5^2$ \\ \hline
$J(6,2)$ & $(15)$ & $A_5^3$ & str. eut. latt. in $\real^{9}$ \\ \hline
$J(7,2)$ & $(21)$ &  $A_6^*$ & str. eut. latt. in $\real^{14}$ \\ \hline
$J(8,2)$ & $(28)$ & $E_7^*$ & str. eut. latt. in $\real^{20}$ \\ \hline
$J(9,2)$ & $(36)$ & $A_8^*$ & str. eut. latt. in $\real^{27}$ \\ \hline
$J(10,2)$ & $(45)$ & $A_9^5$ & str. eut. latt. in $\real^{35}$ \\ \hline
\end{tabular}
\medskip
\caption{Examples of strongly eutactic lattices from Johnson $J(n,2)$ graphs} 
\label{J_graph_ex}
\end{table}
\end{center}
\end{ex}
\vspace{-1cm}
As we we mentioned above, the Johnson graphs $J(n,2)$ are strongly regular, as are their complements Kneser graphs~$KG_{n,2}$. Recall that a (connected) graph $\Gamma$ on $n$ vertices is called {\it strongly regular} with parameters $k$, $\ell$, $m$ whenever it is not complete and:
\begin{enumerate}
\item each vertex is adjacent to $k$ vertices,
\item for each pair of adjacent vertices there are $\ell$ vertices adjacent to both,
\item for each pair of non-adjacent vertices there are $m$ vertices adjacent to both.
\end{enumerate}
Strongly regular graphs are known to have many remarkable properties. In particular, these are precisely the $k$-regular graphs with three distinct eigenvalues. One of these eigenvalues is always $k$ (multiplicity $1$) with the vector $(1,\dots,1)^{\top}$ being a corresponding eigenvector; the other two eigenvalues are roots of the polynomial $x^2 - (\ell-m) x + (m-k)$, which are known to be integers when they have different multiplicity. See Chapter~9 of~\cite{spectra} for many more details.

\begin{ex} \label{other_graphs} We mention a few more examples of notable vertex transitive strongly regular graphs giving rise to interesting lattices (these graphs are described, for instance, in~\cite{spectra} and in~\cite{vt_graph}). These examples are all connected by the common property of being graphs represented by the roots of the lattice~$E_8$ (along with some others already described above; see Section~3.11 of~\cite{distance_regular}, also Section~14.3 of~\cite{deza_book}).

The folded $5$-cube obtained by identifying the antipodal vertices of the $5$-cube is a distance transitive and strongly regular graph on $16$ vertices with parameters $k=5,\ \ell=0,\ m=2$. Its complement (also distance transitive and strongly regular) is called the Clebsch graph. They each have an eigenvalue of multiplicity $5$ ($-3$ and $2$, respectively), and the corresponding lattice is $D_5^*$, the dual of the root lattice $D_5$, where the lattice family $D_n$ is defined as
$$D_n = \left\{ \bx \in \zed^n : \sum_{i=1}^n x_i \equiv 0\ (\md 2) \right\}.$$

The Shrikhande graph can be constructed as Cayley graph of the group $\zed/4\zed \times \zed/4\zed$, taking elements for vertices and connecting two vertices by an edge if and only if their difference is in $\{ \pm (1,0), \pm (0,1), \pm (1,1) \}$. This graph is a vertex transitive, but not distance transitive, and strongly regular graph on $16$ vertices with parameters $k=6,\ \ell=2,\ m=2$. It has an eigenvalue $2$ of multiplicity 6, and the corresponding lattice is~$D_6^+$, which is an example of one of the lattices
$$D_n^+ = D_n \cup \left( \frac{1}{2} \sum_{i=1}^n \be_i + D_n \right),$$
defined for even~$n$. The complement of the Shrikhande graph (also vertex transitive, but not distance transitive, and strongly regular) has eigenvalue $-3$ with multiplicity $6$ and produces the same lattice. Notice that even though the graphs are not distance transitive, the generated lattice is still strongly eutactic.

The Schl\"afli graph is the complement of the intersection graph of the $27$ lines on a cubic surface. It is a distance transitive and strongly regular graph on $27$ vertices with parameters $k=16,\ \ell=10,\ m=8$ and has eigenvalue $4$ of multiplicity $6$.  Its complement (also distance transitive and strongly regular) has eigenvalue $-5$ with multiplicity $6$. Both of these generate the lattice $E_6^*$, the dual of the root lattice~$E_6$. Recall that the lattice $E_8 = D_8^+$, the lattice $E_7$ is the sublattice of $E_8$ with $x_7 = x_8$, and the lattice $E_6$ is the sublattice of $E_8$ with $x_6 = x_7 = x_8$ (see~\cite{CS_book} for more details).

Finally, the Gosset graph (the only one out of these $E_8$-root graphs which is not strongly regular) is a distance transitive graph on $56$ vertices that can be identified with two copies of the set of edges of the complete graph~$K_8$. Then two vertices from the same copy of $K_8$ are connected by an edge if they correspond to disjoint edges of $K_8$, and two vertices from different copies of $K_8$ are connected by an edge if they correspond to edges that meet in a vertex (see~\cite{distance_regular} for more details). The Gosset graph has eigenvalue $9$ of multiplicity~$7$, generating the lattice~$E_7^*$, the dual of~$E_7$.
\end{ex}

The main purpose of all these examples is to demonstrate that this construction of strongly eutactic lattices from distance transitive (and possibly from vertex transitive) graphs appears to produce a wide range of interesting examples already in low dimensions, and hence may be quite useful in higher dimensions too where a classification of strongly eutactic lattices is not yet available.
\smallskip

We also observe here an interesting connection between contact polytopes of some lattices and graphs generating them. For a lattice $\Lambda$, its contact polytope $C(\Lambda)$ is defined as the convex hull of the set of minimal vectors. The significance of the contact polytope is that its vertices are points on the sphere centered at the origin in the sphere packing associated to $\Lambda$ at which neighboring spheres touch it. Hence the number of vertices of $C(\Lambda)$ is the kissing number of $\Lambda$. The skeleton graph of this polytope $\skel(C(\Lambda))$ is the graph consisting of vertices and edges of $C(\Lambda)$. 

Let us consider an example $\Lambda = E_6^*$. The contact polytope of $E_6^*$ has $54$ vertices, split into $27$ $\pm$ pairs: it is a diplo-Schl\"afli polytope (see~\cite{CS-1}). The prefix ``diplo" means double: for a polytope $\Pi$ a diplo-$\Pi$ polytope is a polytope whose vertices are vertices of $\Pi$ and its opposite $-\Pi$. The Schl\"afli polytope, with Coxeter symbol $2_{21}$, has $27$ vertices corresponding to the $27$ lines on a cubic surface~\cite{cox_2-21}. Its skeleton is the Schl\"afli graph $\Gamma$. By Example~\ref{other_graphs} above, $\Gamma$ has an eigenvalue $4$ of multiplicity $6$, and $L_{\Gamma,4} = E_6^*$.

Here is another example of this dual correspondence. For $\Lambda=E_7^*$, its contact polytope is the Gosset polytope (also called Hess polytope) $3_{21}$, which has $56$ vertices (see~\cite{gosset}, \cite{cox_book}). Its skeleton is the Gosset graph $\Gamma$. As we know from Example~\ref{other_graphs} above, $\Gamma$ has an eigenvalue $9$ of multiplicity $7$, and $L_{\Gamma,9} = E_7^*$.

This kind of correspondence certainly does not work for all strongly eutactic lattices. For instance, the contact polytope of $A_n^*$ is a diplo-simplex (see~\cite{CS-1}), and the skeleton graph of a regular simplex on $n+1$ vertices is the complete graph~$K_{n+1}$. By Lemma~\ref{complete}, $K_{n+1}$ generates $A_n$, but not $A_n^*$. On the other hand, the diplo-simplex for $A_3^*$ is a cube, whose skeleton graph $Q_3$ is isomorphic to $H(3,2)$ and the lattice corresponding to eigenvalue $1$ (or $-1$) is $A_3^*$ (see Example~\ref{hamming} above). It would be interesting to understand this correspondence better.

\bigskip

\section{On the coherence of a lattice}
\label{coh}

We conclude with some remarks on the coherence of lattices and frames and their use in the application of compressed sensing. While this discussion is speculative, we hope it will also draw interesting connections and spark interesting future directions. We start with some definitions. Let $L \subset \real^n$ be a lattice. As usual, let $S(L)$ be the set of minimal vectors of $L$, which come in $\pm$ pairs, and let us write $S^*(L)$ for the subset of $S(L)$ where only one vector of each pair is included. Then any two vectors $\bx,\bwy \in S^*(L)$ are linearly independent, so the angle $\theta(\bx,\bwy)$ between them is in the interval $[\pi/3,2\pi/3]$. Define the {\it coherence} of $L$ to be
$$C(L) := \max \{ |\cos \theta(\bx,\bwy)| : \bx \neq \bwy \in S^*(L) \},$$
then $0 \leq C(L) \leq \frac{1}{2}$. In fact, we can speculate a little more about $C(L)$.

The packing density of $L$ is
$$\delta(L) = \frac{\omega_n |L|^n}{2^n \det(L)},$$
where $\omega_n$ is the volume of a unit ball in~$\real^n$. Suppose that a lattice $L \subset \real^n$ has a basis consisting of minimal vectors $\bb_1 \dots \bb_n \in S(L)$ and let
$$B = \begin{pmatrix} \bb_1 \dots \bb_n \end{pmatrix}$$
be the corresponding basis matrix, then $\|\bb_i\| = |L|$ for each $i$. For each $1 \leq i \leq n-1$, let us write $\nu_i$ as the angle between $\bb_{i+1}$ and the subspace spanned by $\bb_1,\dots,\bb_i$. Then
$$\det(L) = | \det(B) | = |L|^n \prod_{i=1}^{n-1} |\sin \nu_i|,$$
and so
$$\delta(L) = \frac{\omega_n}{2^n \prod_{i=1}^{n-1} |\sin \nu_i|}.$$
Therefore
$$\frac{\omega_n}{2^n \delta(L)} =  \prod_{i=1}^{n-1} |\sin \nu_i| \leq \min_{1 \leq i \leq n-1} | \sin \nu_i |,$$
meaning that
\begin{equation}
\label{cos_bnd}
\max_{1 \leq i \leq n-1} | \cos \nu_i | \leq \sqrt{1 - \left( \frac{\omega_n}{2^n \delta(L)} \right)^2}.
\end{equation}

Now, the larger $\prod_{i=1}^n |\sin \nu_i|$ is, the smaller is $\delta(L)$, and it is known that for a lattice $L$ to be a local minimum of~$\delta$ it has to be (weakly) eutactic, but not perfect (Section~9.4 of~\cite{martinet}). Hence it is natural to expect that $\prod_{i=1}^n |\sin \nu_i|$ will be large on non-perfect eutactic lattices (at least some of the time), meaning that the angles $\nu_1,\dots,\nu_{n-1}$ will be large. This suggests that minimal basis vectors might be closer to orthogonal, and so the coherence of the set of minimal vectors, although possibly large relative to the number of minimal vectors, might be expected to be small. We now demonstrate a couple of non-perfect strongly eutactic lattices with coherence $< 1/2$, which come out of our construction of lattices from graphs.

\begin{ex} \label{3-dim} There are three strongly eutactic lattices in~$\real^3$ (up to similarity): $\zed^3 = L_{0_3,0}$, $A_3 = L_{K_3,-1}$ and $A_3^* = L_{H(3,2),1}$, out of which $A_3$ is the only one that is perfect, and hence a local maximum of the packing density function $\delta$ on the space of lattices. Then $\zed^3$ and $A_3^*$ are local minima of~$\delta$. Notice that $S^*(\zed^3)$ is an orthogonal basis, while $S^*(A_3)$ and $S^*(A_3^*)$ are tight frames of cardinalities $6$ and $4$, respectively. The lattice $A_3^*$ can be represented in~$\real^3$ as
$$\begin{pmatrix} \ \ 1 & -1 & \ \ 1 \\ -1 & \ \ 1 & \ \ 1 \\ -1 & -1 & \ \ 1 \end{pmatrix} \zed^3$$
with the set of minimal vectors $\{ (\pm 1, \pm 1, \pm 1) \}$. Hence the coherence~$C(A_3^*) = 1/3$. On the other hand,~$C(A_3) = 1/2$, and no subset of~$S^*(A_3)$ of cardinality $4$ has lower coherence. More generally, the lattice $A_k^*$ represented in~$\real^{k+1}$ has $S^*(A_k^*)$ given by~\eqref{cyclic}, i.e. gives a cyclic $(k,k+1)$-ETF with coherence~$1/k$ discussed in~\cite{etf1}.

Let us also consider the non-perfect strongly eutactic lattice $D_6^+$, generated by the Shrikhande graph. It has $32$ minimal vectors of the form 
$$\frac{1}{2} (\pm 1, \pm 1, \pm 1, \pm 1, \pm 1, \pm 1)$$
with an even number of negative coordinates (\cite{martinet}, Section~4.4), hence $|D_6^+| = \sqrt{3/2}$ and $C(D_6^+) = 1/3$. Thus $S^*(D_6^+)$ is a tight frame of $16$ vectors in~$\real^6$ with coherence~$1/3$: this, again, is an ETF discussed in~\cite{etf1}. This lattice also has a basis of minimal vectors by Theorem~1.1 of~\cite{martinet-bases2}.

Three other examples constructed in~\cite{etf1} we briefly mention are strongly eutactic non-perfect lattices in dimensions $5$, $13$ and $25$, generated by $(10,5)$, $(26,13)$ and $(50,25)$ ETFs, respectively. In all of these three cases the set of minimal vectors of the resulting lattice consists precisely of $\pm$ vectors of the generating frame, and the resulting coherences of these lattices are $1/3$, $1/5$ and $1/7$, respectively. For comparison, the coherence of densest known lattices in dimensions $5$ and $13$ is $1/2$.
\end{ex}

It would be interesting to further investigate coherence of eutactic lattices and, more generally, well-rounded lattices: recall that a full-rank lattice in~$\real^n$ is called well-rounded if it has~$n$ linearly independent minimal vectors; all eutactic and perfect lattices are well-rounded. 
\smallskip

Coherence plays an important role in many applications, and lattice generating ETFs with small coherence are particularly useful. For example, the field of \textit{compressed sensing} aims to recover a sparse vector from a small number of linear measurements. The applications are abundant, ranging from medical imaging and environmental sensing to radar and communications ~\cite{RefWorks:45,RefWorks:373}. Here, we say a vector is $s$-sparse when it has at most $s$ non-zero entries. Put succinctly, compressed sensing aims to recover an $s$-sparse vector $x\in\mathbb{R}^n$  from the measurements $y=Ax\in\mathbb{R}^k$, where $A$ is a suitable $k\times n$ measurement matrix. It is now well known that an $s$-sparse vector $x$ can be efficiently and robustly recovered from measurements $y$ when the number of measurements $k$ is approximately $s\log n$, yielding a significant reduction in the dimension of the representation from $n$ to $s\log n$ (since $s$ is typically much smaller than $n$). For such techniques, one typically constructs $A$ randomly and/or asks that the matrix has highly incoherent columns; this is equivalent to requiring $C(L)$ to be small in situations when columns of $A$ are minimal vectors of a lattice $L$. To this end, it is very natural to consider ETFs and other frames with nice algebraic properties as suitable measurement operators ~\cite{tsiligianni2014construction,FukNedSud18}. Moreover, in many applications, more is known about the signal than simple sparsity; for example, the signal may often also have integer-valued entries or entries in some other lattice.  Such is the case for example in wireless communications ~\cite{rossi2014spatial}, collaborative filtering ~\cite{davenport2016overview}, error correcting codes ~\cite{candes2005error}, and many others.  Although there is some preliminary work for this setting ~\cite{mangasarian2011probability,RefWorks:288,stojnic2010recovery,tian2009detection,zhu2011exploiting,flinth2018promp}, there is still not a rigorous understanding of when and how the lattice structure of the signal can actually be utilized in reconstruction. 

Our work may shed some light on integer-valued sparse recovery by observing the following. If the integer span of an
ETF or another suitable frame is a lattice, then viewing this frame as a measurement matrix (whose columns are the frame vectors), its image restricted to integer-valued signals forms a lattice.  This allows for separation of such images of sparse signals, analogous to the
well-known Johnson-Lindenstrauss lemma, which has been used to guarantee accurate recovery in compressed sensing ~\cite{RefWorks:86}. In fact, when the minimal vectors of the lattice contain the frame vectors, this separation can be bounded.  Viewed in this context, Theorem \ref{rational} gives an answer as to which measurement matrices (given as tight frames) map integer-valued signals to elements of a lattice. Studies of properties of such lattices (e.g. Voronoi cell) have the potential to give stronger guarantees in the integer sparse regime for reconstruction. Of course the integer span of vectors is a larger subset than the image of sparse vectors, however it may be interesting future work to specialize these questions to integer vectors that are in particular also sparse. Group frames may also be interesting for further study given the advantage they give due to their compact representation: fixing a group and picking a starting vector, the entire frame can be generated as its orbit under the group action.

\begin{figure}
\includegraphics[width=2in]{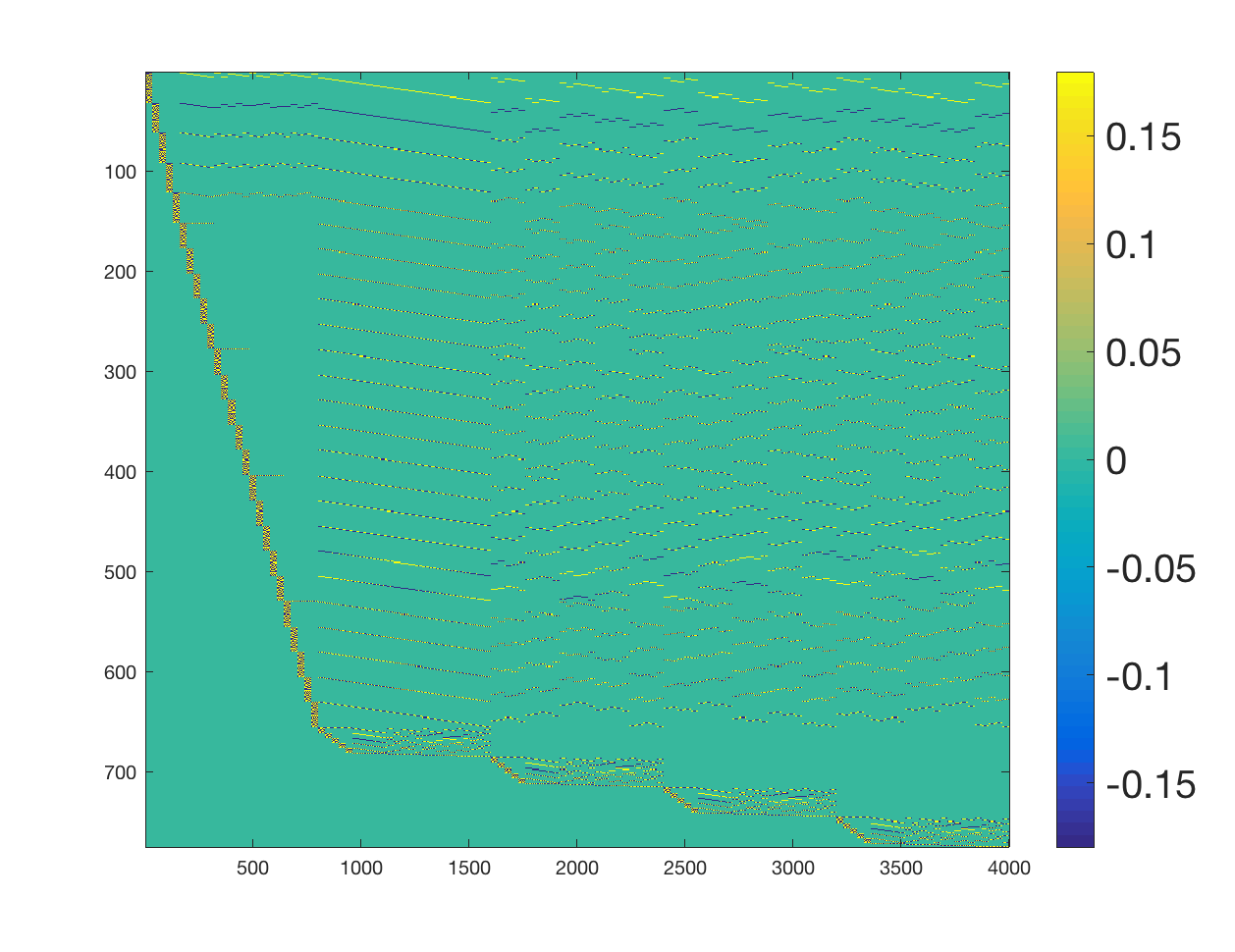}%$\quad$
\includegraphics[width=2in]{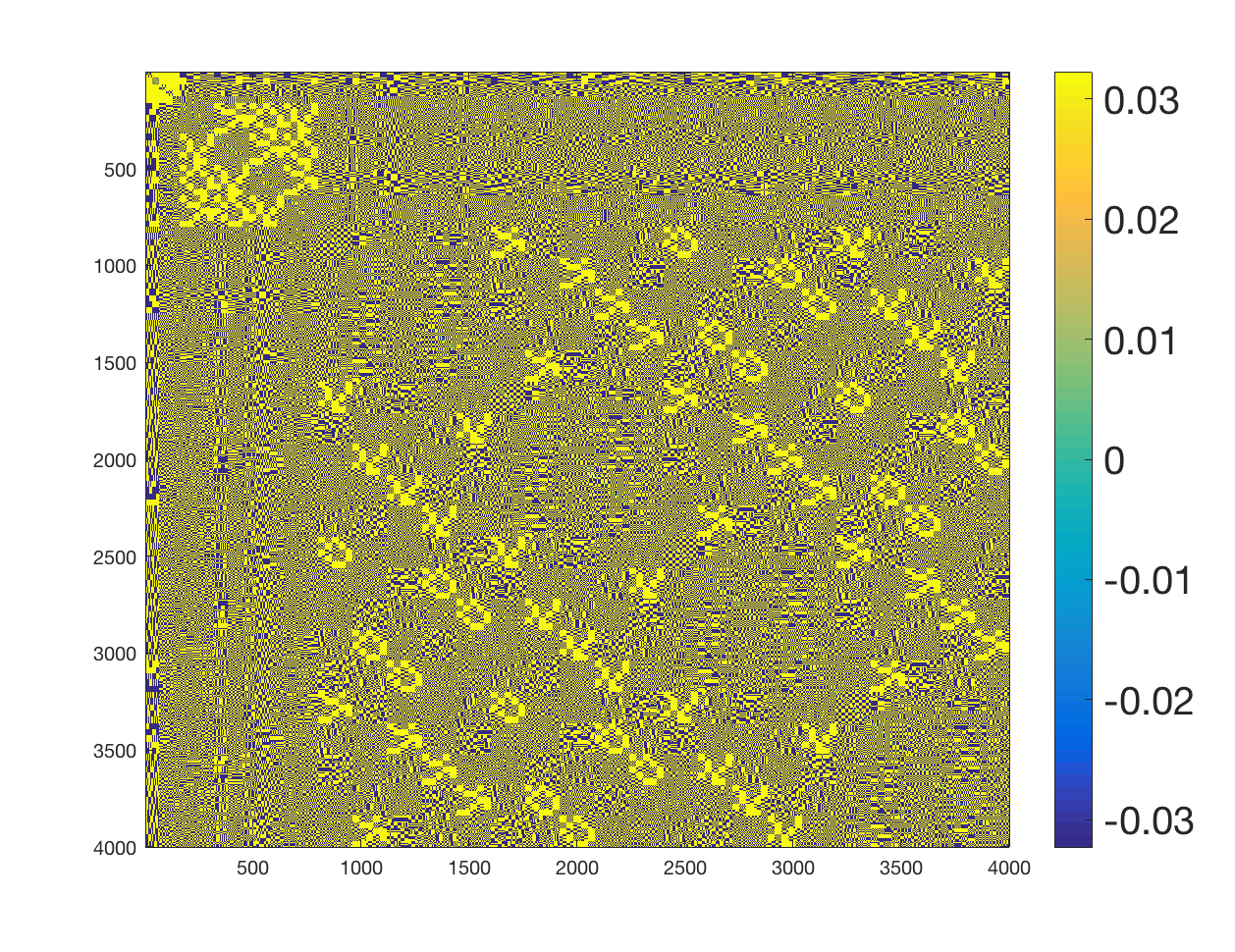}
\caption{Left: A plot of entries in the Steiner ETF.  Right: The corresponding `hollow' Gram matrix ($A^{\top}A-I$).  }\label{theetf}
\end{figure}

To examine how deterministic low-coherence measurement matrices perform in the integer sparse framework we perform a simple experiment using a Steiner ETF of $4000$ vectors in $\real^{775}$, generated from the incidence matrix of an affine Steiner triple system. A schematic representation of this ETF and its Gram matrix is shown in Figure~\ref{theetf}. We chose this measurement matrix for these experiments for a couple of reasons. Steiner ETFs, ETFs generated from a type of combinatorial construction, have been singled out as some of the ETFs with the most potential in application to problems in compressed sensing \cite{fickus2012steiner}. These Steiner ETFs stand out because by working in a sufficiently large dimension the coherence can be made arbitrarily small and the redundancy as large as desired, this property being inherited from known constructions of Hadamard matrices and Steiner triple systems used to generate these incoherent frames \cite{fickus2012steiner,goethals}. Although these matrices have other undesirable properties such as being sparse themselves, the freedom to generate large matrices with small coherence is instrumental in sparse recovery given the well-studied relation between low-coherence matrices and guarantees in compressed sensing.

\begin{figure}
\includegraphics[width=5in]{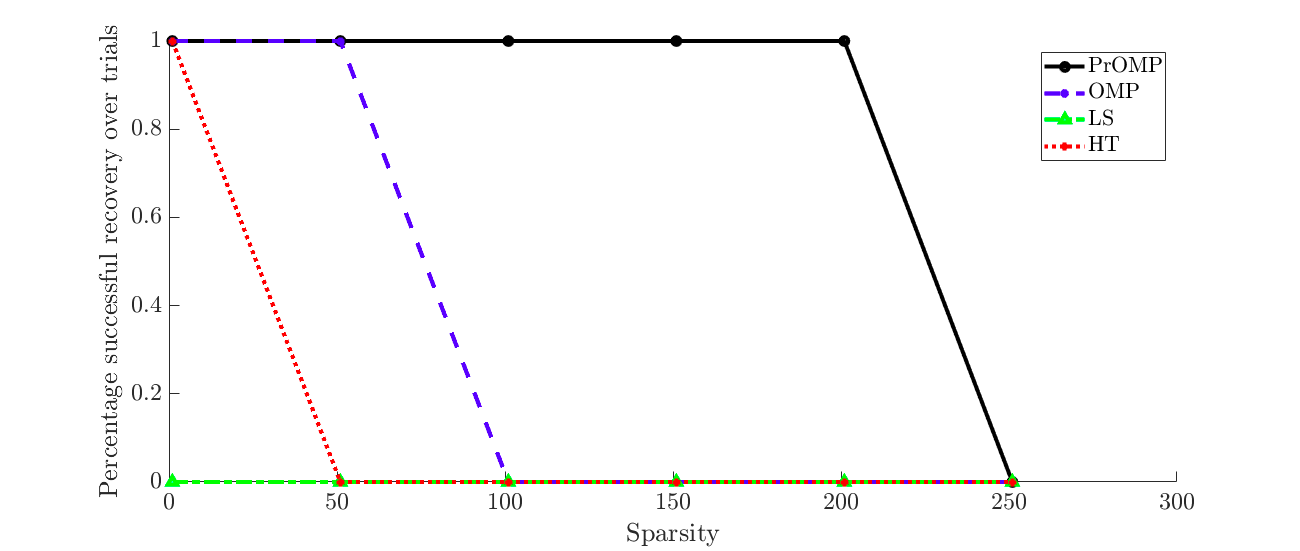} \\ 

\includegraphics[width=5in]{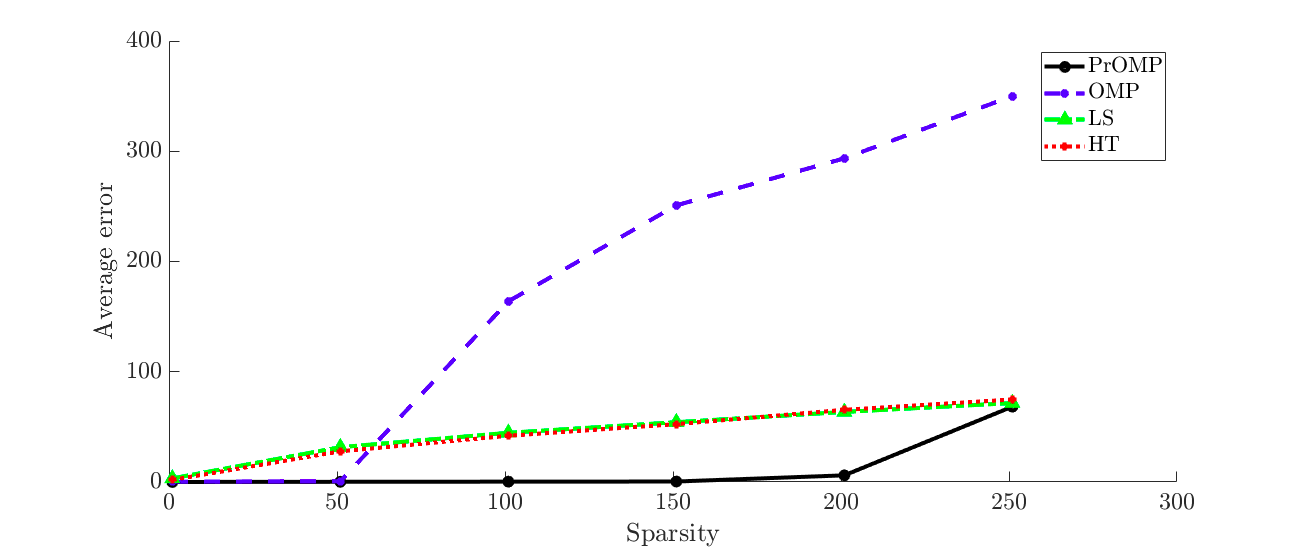}
\caption{Recovery results for various algorithms (PrOMP, OMP, Hard Thresholding, Least Squares) using a Steiner ETF in $\real^{775}$, size $4000$, as the measurement matrix. Left: Percentage of accurate recovery. Right: Noise added to the measurements to have norm~$0.1$.}\label{fig:exp}
\end{figure}

Denoting this frame of vectors by $F$, we acquire the measurements $y = Fx$ or the noisy measurements $y = Fx+e$ where $x$ is a vector of varying sparsity and $e$ is scaled Gaussian noise. We then use various compressed sensing algorithms to recover $\hat{x}$ and calculate how often recovery is exact ($x=\hat{x}$) in the noiseless case, and the magnitude of the recovery error ($\|x-\hat{x}\|_2$) in the noisy case. We show results for the simple least-squares method (LS) that simply sets $\hat{x} = F^\dagger y$, basic hard thresholding (HT) which first estimates the support of $x$ via the proxy $F^T y$ and then performs least-squares over that support, Orthogonal Matching Pursuit \cite{Paper9} (OMP) which is an iterative greedy algorithm, and PrOMP \cite{flinth2018promp} which is a modification of OMP for integer-valued signals.  The results are shown in Figure \ref{fig:exp}, where we see unsurprisingly that PrOMP performs quite well in this case, confirming the previous observations of effectiveness of pre-processing steps in lattice-valued compressed sensing. The previous analysis in \cite{flinth2018promp} has explained via a concentration of measure argument why this should hold for Gaussian matrices, but numerically there is some evidence that performance improvements hold for deterministic measurements and integer signals in iterative compressed sensing procedures when a pre-processing step, as is found in PrOMP, is applied.

\bigskip

\noindent
{\bf Acknowledgement:} We are grateful to Benny Sudakov, Achill Sch\"urmann, and Shahaf Nitzan for helpful discussions. We also thank Mathieu Dutour Sikiri\'c for helping to verify our example computations. Finally, we want to thank the referees for comments and corrections that greatly improved the quality of the paper.

%\nocite{*}
\bibliographystyle{myalpha}  % Here the bibliography 
\bibliography{g_frames}    % is inserted.
\end{document}